\documentclass[11pt]{article}
\usepackage[T1]{fontenc}
\usepackage[utf8]{inputenc}
\usepackage{mathpazo}
\usepackage{microtype}
\usepackage{amsmath,amssymb,amsthm}
\usepackage[margin=1.1in]{geometry}
\usepackage{booktabs}
\usepackage{graphicx}
\usepackage[titles]{tocloft}

\usepackage{xcolor}
\usepackage{sectsty}
\usepackage[colorlinks=true,linkcolor=blue!50!black,citecolor=blue!50!black,urlcolor=blue!50!black]{hyperref}
\allsectionsfont{\normalfont\scshape}
\newtheoremstyle{plainsc}{8pt}{8pt}{\itshape}{}{\scshape}{.}{0.5em}{}
\newtheoremstyle{defsc}{8pt}{8pt}{\normalfont}{}{\scshape}{.}{0.5em}{}
\theoremstyle{plainsc}
\newtheorem{theorem}{Theorem}[section]

\newtheorem{corollary}[theorem]{Corollary}
\theoremstyle{defsc}
\newtheorem{remark}[theorem]{Remark}
\newcommand{\C}{\mathbb{C}}
\newcommand{\R}{\mathbb{R}}
\newcommand{\E}{\mathbb{E}}
\DeclareMathOperator{\Real}{Re}
\DeclareMathOperator{\Imag}{Im}

\title{On Hilbert's 8th Problem}
\author{Nicholas G. Polson\\ \small Booth School of Business, University of Chicago}
\date{October 2026}

\begin{document}
\maketitle

\begin{abstract}
Hilbert's eighth problem asks, among other things, for a proof of the Riemann Hypothesis. We approach it through the probability law on $[0,\infty)$ whose Laplace transform is $\varphi(s) = \xi(1/2)/\xi(1/2+\sqrt s)$. Unconditionally this law is self decomposable, with Lévy density $\Theta(x)/x$ where $\Theta(x) = \sum e^{(\rho-1/2)^2x}$ is the heat trace over the zeros; we write the Lévy measure out in three closed forms, through the zeros, through the primes and the gamma factor, and through Brownian first passage densities. RH is equivalent to the law being a generalized gamma convolution, that is, to the existence of a positive Thorin measure, which is then the counting measure of the squared ordinates. Moving the base point to $\alpha$, the reciprocal $\xi(\alpha)/\xi(\alpha+\sqrt s)$ is a generalized gamma convolution if and only if no zero lies to the right of $\alpha$, unconditionally for $\alpha \ge 1$, with Thorin measure the zero counting measure smoothed by Cauchy kernels, whose density in the ordinate variable is $\pi^{-1}\Real\,\xi'/\xi$ on the line $\Real s = \alpha$. By Bondesson's continuity theorem RH follows if these laws survive on a sequence $\alpha_n \downarrow 1/2$; this survival is a backward Cauchy flow, equivalent at any fixed $\alpha$ to the quasi Riemann hypothesis and, for all $\alpha$, to Lagarias' positivity criterion. We separate throughout what is unconditional from what is equivalent to RH, and identify the one sided bound on $\zeta'/\zeta$ that any proof along these lines must establish. At $\alpha = 1$ this bound holds, and we prove its smoothed form from the primes, via Mertens' theorem and with no information about the zeros beyond the prime number theorem, on every window of one mean zero spacing; passing to finer windows is where the zeros reenter.
\end{abstract}



\section{Introduction}\label{sec:intro}

Hilbert's eighth problem, posed in 1900, groups together several questions on prime numbers, the first and most famous being the Riemann Hypothesis. These notes approach it from probability. The starting point is a classical fact about the completed zeta function $\xi(z) = \tfrac12z(z-1)\pi^{-z/2}\Gamma(z/2)\zeta(z)$: the function $\xi(1/2+w)$ is even and entire of order one, so $F(s) = \xi(1/2+\sqrt s)/\xi(1/2)$ is entire of order one half, and its reciprocal $\varphi = 1/F$ is a candidate Laplace transform. Under RH it is the Laplace transform of $\sum_\gamma E_\gamma/\gamma^2$ with independent standard exponentials $E_\gamma$, a generalized gamma convolution (GGC) in the sense of Thorin and Bondesson \cite{Bondesson}, whose Thorin measure is the counting measure of the squared ordinates $\gamma^2$. The question is how much of this survives without RH, and what the surviving structure says about RH itself.

The answer has three layers. First, unconditionally $\varphi$ is the Laplace transform of a self decomposable law whose canonical function is the heat trace $\Theta(x) = \sum e^{(\rho-1/2)^2x}$ over the zeros (Theorem~\ref{thm:law}); the Lévy measure is written out in Section~\ref{sec:levy} through the zeros, through the explicit formula, and through Brownian first passage densities at the levels $\log n$. Second, RH is equivalent to complete monotonicity of $\Theta$ (Theorem~\ref{thm:rh}), so the Thorin measure of Section~\ref{sec:thorin} exists as a positive measure exactly under RH; the 2 sin$^2$ lemma of Section~\ref{sec:sin2} explains how a Lévy measure on the primal side becomes a Thorin measure on the reciprocal side, and where this transfer breaks for $\xi$. Section~\ref{sec:division} separates, with proofs, everything that is unconditional from everything that is equivalent to RH. Third, moving the base point gives a family of laws $\mu^{(\alpha)}$ with Laplace transforms $\xi(\alpha)/\xi(\alpha+\sqrt s)$ that are GGCs, unconditionally, for $\alpha \ge 1$ (Theorem~\ref{thm:base}); by Bondesson's continuity theorem (Theorem~\ref{thm:continuity}) RH follows if they remain GGCs on a sequence $\alpha_n \downarrow 1/2$. The descent is a backward Cauchy flow on the zero counting measure, and Theorem~\ref{thm:lagarias} identifies its survival with the quasi Riemann hypothesis at each fixed $\alpha$ and with Lagarias' positivity of $\Real\,\xi'/\xi$ for all $\alpha$. Section~\ref{sec:alpha1} treats the base point $\alpha = 1$ on its own: there the Thorin measure is absolutely continuous exactly because $\zeta$ has no zeros on the line $\Real s = 1$, which is the prime number theorem in Thorin form; the laws to the right of $1$ are Cauchy smoothings of the law at $1$ (Theorem~\ref{thm:poisson}), and the positivity of the Thorin density, smoothed at the scale of the mean zero spacing, follows from Mertens' theorem with no information about the zeros beyond the prime number theorem (Theorem~\ref{thm:fejer}). Section~\ref{sec:path} draws the consequences for a proof, and Section~\ref{sec:open} lists open problems.

The emphasis throughout is on what is proved. No step here proves RH, and the final sections say precisely which inequality, a one sided bound on $-\Real\,\zeta'/\zeta$ at some fixed $\alpha < 1$ and all heights, a proof along these lines would require. Proofs not given here, and the numerical scripts, are in the companion notes \cite{Ladder}.

\section{The law and the dictionary}\label{sec:law}

\emph{Notation.} Throughout, $\rho = \beta + i\gamma = \tfrac12 + \delta + i\gamma$ runs over the nontrivial zeros of $\zeta$ with $\gamma > 0$, counted with multiplicity; $N(T)$ counts them up to height $T$ and $N_0(T)$ counts those on the line. We write $\gamma_\rho = \gamma - i\delta$ and $s_\rho = (\rho-\tfrac12)^2 = -\gamma_\rho^2$. The reciprocal Xi law is $\mu$, its Laplace transform $\varphi$, its exponent $\psi = -\log\varphi$, and its canonical function $\Theta$; $T_k = (-1)^k\Theta^{(k)}$. At a base point $\alpha$ the law is $\mu^{(\alpha)}$, with Thorin measure $U_\alpha$ and widths $c_\rho = \alpha - \beta$. $H = 3\cdot10^{12}$ is the verified height, $\gamma_{\rm E}$ is Euler's constant, $\Lambda$ is the von Mangoldt function, and $\Lambda$ without argument is the de Bruijn Newman constant.

Since $\xi(1/2+w)$ is even and entire of order one in $w$, the function $F(s) = \xi(1/2+\sqrt s)/\xi(1/2)$ is entire of order one half in $s$, with $F(0) = 1$, $F > 0$ on $[0,\infty)$, and Hadamard product
\begin{equation}\label{eq:hadamard}
F(s) = \prod_{\rho}\Big(1 - \frac{s}{s_\rho}\Big), \qquad s_\rho = (\rho - \tfrac12)^2,
\end{equation}
one factor per pair $\{\rho, 1-\rho\}$ of nontrivial zeros, equivalently one per zero $\rho = \tfrac12 + \delta + i\gamma$ in the upper half plane. Write $\gamma_\rho = -i(\rho - 1/2) = \gamma - i\delta$, so that $s_\rho = -\gamma_\rho^2$.

\begin{theorem}\label{thm:law}
Unconditionally, $\varphi = 1/F$ is the Laplace transform of a self decomposable probability law $\mu$ on $[0,\infty)$:
\[
\varphi(s) = \frac{\xi(1/2)}{\xi(1/2+\sqrt s)} = \int_0^\infty e^{-sy}\,\mu(dy) = \exp\Big(-\int_0^\infty(1-e^{-sx})\,\frac{\Theta(x)}{x}\,dx\Big), \qquad s \ge 0,
\]
with no drift and no Gaussian part, where the canonical function
\begin{equation}\label{eq:theta}
\Theta(x) = \sum_{\gamma>0}e^{-\gamma_\rho^2 x} = \sum_{\gamma>0}e^{(\rho-1/2)^2x}
\end{equation}
is positive and nonincreasing on $(0,\infty)$. The law $\mu$ is absolutely continuous and unimodal. It is a generalized gamma convolution if and only if RH holds, and then
\[
\mu = \text{law of } \sum_{\gamma>0}\frac{E_\gamma}{\gamma^2}, \qquad E_\gamma \text{ independent standard exponentials.}
\]
\end{theorem}

The three claims about $\Theta$ are: $\Theta > 0$ and $-\Theta' > 0$ unconditionally, which give infinite divisibility and self decomposability; and $\Theta$ completely monotone if and only if RH, which gives the GGC property. The first two are proved in Section~\ref{sec:division} from the verified zeros and an explicit zero counting bound, following \cite{Ladder}; the third is proved now, after the dictionary that translates properties of the law into properties of $\Theta$.

The whole subject rests on one identity. Let $\psi = \log F = -\log\varphi$. For $\Real s > -\gamma_1^2$,
\begin{equation}\label{eq:dict}
\psi'(s) = \frac{F'}{F}(s) = \sum_\rho\frac{1}{s - s_\rho} = \int_0^\infty e^{-sx}\,\Theta(x)\,dx,
\end{equation}
by termwise integration of $e^{-sx}e^{s_\rho x}$, justified since $\Theta(x) = O(x^{-1/2}\log(1/x))$ near $0$ and decays like $e^{-\gamma_1^2x}$. Integrating from $0$, $\psi(s) = \int_0^\infty(1-e^{-sx})\Theta(x)x^{-1}dx$. Three classical facts then translate properties of the law into properties of $\Theta$ \cite{Bondesson,SSV}:
\begin{enumerate}
\item[(a)] $\varphi$ is the Laplace transform of an infinitely divisible law on $[0,\infty)$ iff $\psi$ is a Bernstein function iff $\psi'$ is completely monotone iff $\Theta \ge 0$, by Bernstein's theorem applied to \eqref{eq:dict} and uniqueness of Laplace transforms.
\item[(b)] The law is self decomposable iff its Lévy density has the form $k(x)/x$ with $k$ nonincreasing, iff $-\Theta' \ge 0$.
\item[(c)] The law is a GGC iff $k = \Theta$ is completely monotone, iff $\psi'$ is a Stieltjes function, and then the Thorin measure $U$ is the representing measure $\Theta(x) = \int e^{-tx}U(dt)$.
\end{enumerate}
With $T_k(x) = (-1)^k\Theta^{(k)}(x)$, every statement in these notes is a statement about the sign of $T_k$ on some range of $(k,x)$: (a) is $T_0 \ge 0$, (b) is $T_1 \ge 0$, (c) is $T_k \ge 0$ for all $k$. The next theorem is the top of that ladder.

\begin{theorem}\label{thm:rh}
RH $\iff$ $\Theta$ is completely monotone $\iff$ $\varphi$ is a GGC Laplace transform.
\end{theorem}
\begin{proof}
Under RH, $s_\rho = -\gamma^2$, so $\Theta(x) = \sum e^{-\gamma^2x} = \int e^{-tx}\,dU(t)$ with $U = \sum\delta_{\gamma^2}$, completely monotone. Conversely, if $\Theta$ is completely monotone then $\Theta = \int e^{-tx}\mu(dt)$ with $\mu \ge 0$ by Bernstein, so by \eqref{eq:dict} $\psi'(s) = \int\mu(dt)/(t+s)$ is a Stieltjes function and extends holomorphically to $\C\setminus(-\infty,0]$. But $\psi' = F'/F$ has a pole at every $s_\rho$. Hence every $s_\rho = (\rho-\tfrac12)^2$ is a nonpositive real number, so $\rho - \tfrac12$ is purely imaginary (it is not real, since $\xi$ has no real zeros), which is RH.
\end{proof}
The proof shows exactly what RH buys: it places every pole of $\psi'$ on the cut, and only then can $\psi'$ be a Stieltjes function. Without RH the poles at $s_\rho$, $\overline{s_\rho}$ off the cut are exactly what Stieltjes functions cannot have.

\section{The Lévy measure}\label{sec:levy}

The Lévy measure of $\mu$ is $\nu(dx) = \Theta(x)\,x^{-1}dx$, so describing $\nu$ means describing $\Theta$. This section gives three closed forms of $\Theta$. Each exposes a different structure: the first the zeros, the second the primes and the gamma factor, the third a continuum of Brownian levels. It then records the behaviour of $\Theta$ at the two ends of $(0,\infty)$ and the moments of the law.

\subsection*{First form: the zeros}

This is the definition \eqref{eq:theta}. A zero $\tfrac12 + \delta + i\gamma$ contributes $e^{(\delta^2-\gamma^2)x}\cos(2\delta\gamma x)$, and its partner $\tfrac12 - \delta + i\gamma$ contributes the conjugate, so $\Theta$ is real. Under RH every $\delta$ vanishes and
\begin{equation}\label{eq:zeroform}
\Theta(x) = \sum_{\gamma>0}e^{-\gamma^2x}, \qquad \nu(dx) = \sum_{\gamma>0}\frac{e^{-\gamma^2x}}{x}\,dx .
\end{equation}
Each term is the Lévy density of one exponential variable $E_\gamma/\gamma^2$, which is how the GGC representation $\sum E_\gamma/\gamma^2$ appears.

\subsection*{Second form: the explicit formula}

The Weil explicit formula~\cite{Weil} with the Gaussian test function $e^{-r^2x}$ expresses the same sum through primes. It is an identity, with no hypothesis:
\begin{equation}\label{eq:weil}
\Theta(x) = e^{x/4} - \frac{\log\pi}{4\sqrt{\pi x}} + \frac{1}{4\pi}\int_\R e^{-r^2x}\,\Real\frac{\Gamma'}{\Gamma}\Big(\frac14+\frac{ir}{2}\Big)dr - \frac{1}{2\sqrt{\pi x}}\sum_{n\ge2}\frac{\Lambda(n)}{\sqrt n}\,e^{-(\log n)^2/(4x)} .
\end{equation}
Read from left to right, the terms are the pole of $\zeta$, the archimedean place (two terms), and the finite places. Numerically the two forms agree to every digit computed; for example both give $\Theta(0.003) = 1.09522012145284\ldots$.

\subsection*{Third form: Brownian levels}

The prime sum in \eqref{eq:weil} is a sum over the levels $\log n$. The third form rewrites the archimedean integral in the same language, as a continuum of levels $y > 0$ with density
\[
m(y) = \frac{e^{-y/2}}{1-e^{-2y}}, \qquad\text{and constant}\qquad \kappa = \tfrac12\Big(\log\pi - \frac{\Gamma'}{\Gamma}\big(\tfrac14\big)\Big) = 2.686091\ldots .
\]
Then
\begin{equation}\label{eq:levelform}
\Theta(x) = e^{x/4} - \frac{\kappa}{2\sqrt{\pi x}} + \frac{1}{2\sqrt{\pi x}}\int_0^\infty m(y)\Big(1-e^{-y^2/(4x)}\Big)dy - \frac{1}{2\sqrt{\pi x}}\sum_{n\ge2}\frac{\Lambda(n)}{\sqrt n}e^{-(\log n)^2/(4x)} .
\end{equation}
To derive it, insert Gauss's integral $\frac{\Gamma'}{\Gamma}(z) = \int_0^\infty\big(e^{-u}u^{-1} - e^{-zu}(1-e^{-u})^{-1}\big)du$ into \eqref{eq:weil}, integrate in $r$ with $\int_\R e^{-r^2x}\cos(ru/2)\,dr = \sqrt{\pi/x}\,e^{-u^2/(16x)}$, and substitute $u = 2y$; the constant $\kappa$ collects $\frac{\Gamma'}{\Gamma}(\tfrac14) - \log\pi = -2\kappa$.

Dividing by $x$ turns \eqref{eq:levelform} into a statement about Brownian motion. Write
\[
f_a(x) = \frac{a\,e^{-a^2/(4x)}}{2\sqrt\pi\,x^{3/2}}
\]
for the density of the first passage time of a Brownian motion with variance $2$ to the level $a$. Then
\begin{equation}\label{eq:levydensity}
\frac{\Theta(x)}{x} = \frac{e^{x/4}}{x} \;-\; \frac{\kappa}{2\sqrt\pi\,x^{3/2}} \;+\; \int_0^\infty\frac{m(y)}{y}\Big(\frac{y}{2\sqrt\pi\,x^{3/2}} - f_y(x)\Big)dy \;-\; \sum_{n\ge2}\frac{\Lambda(n)}{\sqrt n\,\log n}\,f_{\log n}(x) .
\end{equation}
The four pieces have the following meanings.
\begin{itemize}
\item[] \emph{The pole.} $e^{x/4}/x$, the Lévy density of an exponential law with negative rate. Its Laplace exponent $\int(1-e^{-sx})e^{x/4}x^{-1}dx$ is formally $\log(1-4s)$, which is the factor $z(z-1) = s - \tfrac14$ of $\xi$ at $z = \tfrac12+\sqrt s$, up to a constant.
\item[] \emph{The stable correction.} minus the Lévy density of a positive stable law of index one half, with Laplace exponent $-\kappa\sqrt s$.
\item[] \emph{The archimedean levels.} for each level $y$, a stable density minus a first passage density, which is a positive function; the integral against $m(y)/y$ has Laplace exponent $\log\Gamma(\tfrac14+\tfrac{\sqrt s}{2}) - \log\Gamma(\tfrac14) - \tfrac{\sqrt s}{2}\frac{\Gamma'}{\Gamma}(\tfrac14)$.
\item[] \emph{The primes.} minus a superposition of first passage densities at the levels $\log n$, with weights $\Lambda(n)/(\sqrt n\log n)$; this is the Euler product in Lévy form, with Laplace exponent formally $\log\zeta(\tfrac12+\sqrt s) - \log\zeta(\tfrac12)$.
\end{itemize}
The first and last exponents are formal, because $\int_1^\infty e^{x/4}x^{-1}dx$ and $\sum\Lambda(n)/(\sqrt n\log n)$ both diverge. The divergences cancel against each other, and the four exponents add up to the entire function $\log F(s) = \log\xi(\tfrac12+\sqrt s) - \log\xi(\tfrac12)$, as they must.

\subsection*{The two ends and the moments}

As $x \to 0^+$ the prime sum in \eqref{eq:weil} is exponentially small and the archimedean integral dominates:
\begin{equation}\label{eq:small}
\begin{gathered}
\Theta(x) = \frac{1}{8\sqrt{\pi x}}\log\frac{1}{16\pi^2e^{\gamma_{\rm E}}x} + 1 + C_0 + o(1),\\
C_0 = \frac{1}{4\pi}\int_\R\Big(\Real\frac{\Gamma'}{\Gamma}\Big(\frac14+\frac{ir}2\Big) - \log\frac{|r|}{2}\Big)dr = -0.1227 .
\end{gathered}
\end{equation}
So near the origin the Lévy density behaves like $x^{-3/2}\log(1/x)/(8\sqrt\pi)$: a positive stable law of index one half with a logarithmic correction. As $x\to\infty$ the first zero dominates, $\Theta(x) \sim e^{-\gamma_1^2x}$ with $\gamma_1 = 14.1347\ldots$, so the law has exponential moments up to order $\gamma_1^2 = 199.79$.

\begin{figure}[htbp]
\centering
\includegraphics[width=0.9\textwidth]{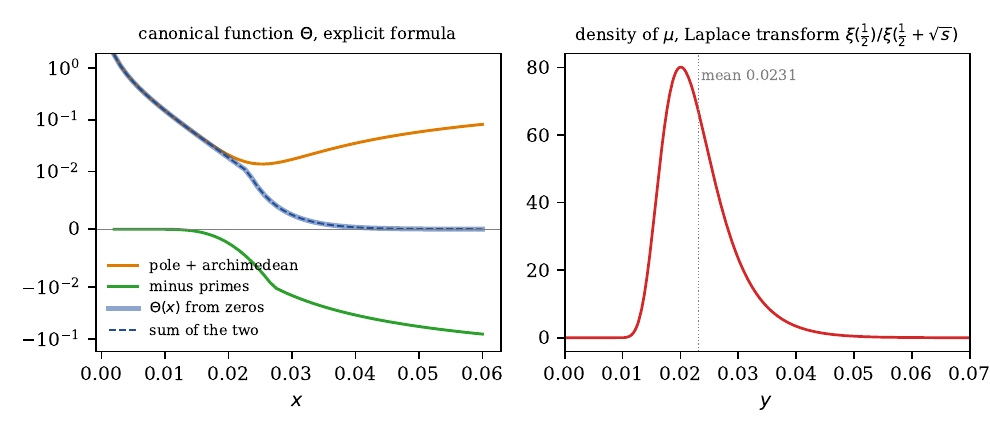}
\caption{Left: the canonical function $\Theta$ as the sum of the pole and archimedean terms (orange) and minus the prime term (green) of the explicit formula, against the zero sum (blue); the vertical scale is logarithmic away from zero. Right: the density of the reciprocal Xi law $\mu$, with its mean $\sum\gamma^{-2} = 0.0231$.}
\label{fig:law}
\end{figure}

The cumulants follow from \eqref{eq:dict}. Since $\log\E\,e^{-sY} = -\log F(s) = \sum_n\kappa_n(-s)^n/n!$,
\[
\kappa_n(\mu) = (n-1)!\int_0^\infty x^{n-1}\Theta(x)\,dx = (n-1)!\sum_{\gamma>0}\gamma_\rho^{-2n} = (-1)^{n+1}n!\,a_n ,
\]
where $a_n$ is the coefficient of $s^n$ in $\log F(s)$; numerically $a_1 = 0.0231049931$, $a_2 = -1.85863\times10^{-5}$, $a_3 = 4.80580\times10^{-8}$ and $a_4 = -1.65758\times10^{-10}$. Hence
\begin{gather*}
\kappa_1 = \sum\gamma^{-2} = \frac{\xi''(1/2)}{2\xi(1/2)} = 0.0231049931,\qquad \kappa_2 = \sum\gamma^{-4} = 3.71726\times10^{-5},\\
\kappa_3 = 2\sum\gamma^{-6} = 2.88348\times10^{-7},\qquad \kappa_4 = 6\sum\gamma^{-8} = 3.97819\times10^{-9},
\end{gather*}
and the sums over the first $600$ zeros, with a Riemann von Mangoldt tail, reproduce these to six digits. The law is tightly concentrated: its mean is $0.0231$, its standard deviation $0.0061$, and its skewness $\kappa_3/\kappa_2^{3/2} = 1.27$. Figure~\ref{fig:law} shows both sides of the picture. On the left, the explicit formula \eqref{eq:weil} at work: the pole and archimedean terms grow with $x$, the prime term switches on near $x \approx 0.02$, when $e^{-(\log 2)^2/(4x)}$ ceases to be negligible, and the two cancel to leave the rapidly decaying heat trace; the sum agrees with the zero side to machine precision. On the right, the density of $\mu$ itself, obtained by Fourier inversion of $t \mapsto \xi(\tfrac12)/\xi(\tfrac12+\sqrt{-it})$ on $0 \le t \le 3\times10^4$, where the transform has decayed to $10^{-45}$. The inversion recovers total mass $0.99997$ and mean $0.023103$, against the exact $0.023105$; the density is unimodal, as self decomposability requires, with mode near $0.020$ and the right skew of a law whose Lévy density is concentrated near the origin.

\begin{table}[htbp]
\centering
\small
\begin{tabular}{rrrr}
\toprule
$x$ & $\Theta(x)$ & Lévy density $\Theta(x)/x$ & driver density $-\Theta'(x)$\\
\midrule
$0.001$ & $3.70458$ & $3704.58$ & $3644.15$\\
$0.003$ & $1.09522$ & $365.073$ & $465.26$\\
$0.01$ & $0.149698$ & $14.9698$ & $33.7292$\\
$0.03$ & $0.00249614$ & $0.0832046$ & $0.499130$\\
$0.1$ & $2.1048\times10^{-9}$ & $2.1048\times10^{-8}$ & $4.2052\times10^{-7}$\\
\bottomrule
\end{tabular}
\caption{The canonical function, the Lévy density and the jump density of the Ornstein Uhlenbeck driver, from the zero form with $600$ zeros.}
\end{table}

Self decomposability has a dynamic reading which belongs here, since it concerns the Lévy side. It means $\mu$ is the stationary law of $dY_t = -Y_t\,dt + dD_t$, equivalently $Y \stackrel{d}{=} \int_0^\infty e^{-t}dD_t$, for a subordinator $D$ whose Lévy measure $\rho$ has tail $\rho((x,\infty)) = \Theta(x)$. So the driver has Lévy density
\begin{equation}\label{eq:driver}
-\Theta'(x) = \sum_{\gamma>0}\gamma_\rho^2\,e^{-\gamma_\rho^2x} \;\stackrel{\rm RH}{=}\; \sum_{\gamma>0}\gamma^2e^{-\gamma^2x},
\end{equation}
infinite activity, total mass $\Theta(0+) = \infty$, and from \eqref{eq:weil} the explicit form
\begin{align*}
-\Theta'(x) &= -\frac{e^{x/4}}{4} - \frac{\log\pi}{8\sqrt\pi\,x^{3/2}} + \frac{1}{4\pi}\int_\R r^2e^{-r^2x}\,\Real\frac{\Gamma'}{\Gamma}\Big(\frac14+\frac{ir}2\Big)dr\\
&\quad + \frac{1}{2\sqrt{\pi x}}\sum_{n\ge2}\frac{\Lambda(n)}{\sqrt n}e^{-(\log n)^2/(4x)}\Big(\frac{(\log n)^2}{4x^2} - \frac{1}{2x}\Big) .
\end{align*}
Under RH, $D$ is the independent superposition over the zeros of compound Poisson processes of rate one with $\mathrm{Exp}(\gamma^2)$ jumps, and the Ornstein Uhlenbeck equilibrium of each is $\mathrm{Exp}(\gamma^2)$, recovering $\sum E_\gamma/\gamma^2$.

\section{The Thorin measure}\label{sec:thorin}

The Lévy measure always exists; the Thorin measure is the finer object whose existence as a positive measure is RH. Recall that a law on $[0,\infty)$ is a GGC when its Laplace exponent has the form
\[
\psi(s) = -\log\varphi(s) = \log F(s) = as + \int_{(0,\infty)}\log\Big(1+\frac{s}{t}\Big)U(dt),
\]
with drift $a\ge0$ and a nonnegative Thorin measure $U$, or equivalently when the canonical function is a Laplace transform $\Theta(x) = \int e^{-tx}U(dt)$. Whether or not $U$ is positive, it can be read off from $F$ by Stieltjes inversion of $\psi' = F'/F$:
\begin{equation}\label{eq:inversion}
U(dt) = \lim_{\varepsilon\to0^+}-\frac1\pi\Imag\frac{F'}{F}(-t+i\varepsilon)\,dt, \qquad U\big((0,t]\big) = \frac1\pi\arg F(-t+i0),
\end{equation}
where the argument is followed continuously from $s=0$ along the upper edge of the negative axis. In words: the Thorin distribution function is the argument of $F$ on the cut, divided by $\pi$.

\subsection*{Under RH}

Then everything is explicit:
\begin{equation}\label{eq:thorinRH}
U = \sum_{\gamma>0}\delta_{\gamma^2}, \qquad a = 0, \qquad U\big((0,t]\big) = N(\sqrt t),
\end{equation}
where $N(T)$ counts the zeros with $0 < \gamma \le T$. Table~\ref{tab:count} checks this through \eqref{eq:inversion}, evaluated at $\varepsilon = 1$ rather than $0^+$ so that the atoms are slightly smoothed; the agreement is to the expected $10^{-2}$.

\begin{table}[htbp]
\centering
\small
\begin{tabular}{rrrrrrrrr}
\toprule
$\sqrt t$ & 15 & 22 & 26 & 31 & 35 & 40 & 45 & 50\\
\midrule
$\pi^{-1}\arg F(-t+i)$ & 0.995 & 1.999 & 2.998 & 3.997 & 5.002 & 6.006 & 8.001 & 9.988\\
$N(\sqrt t)$ & 1 & 2 & 3 & 4 & 5 & 6 & 8 & 10\\
\bottomrule
\end{tabular}
\caption{The Thorin distribution function as a boundary argument, against the zero count.}
\label{tab:count}
\end{table}

The Riemann von Mangoldt formula $N(T) = \vartheta(T)/\pi + 1 + S(T)$, with $\vartheta$ the Riemann Siegel theta function and $S(T) = \pi^{-1}\arg\zeta(\tfrac12+iT)$, splits $U$ into a smooth part and a fluctuation:
\begin{equation}\label{eq:thorinsmooth}
U\big((0,t]\big) = \frac{\vartheta(\sqrt t)}{\pi} + 1 + S(\sqrt t), \qquad U(dt) = \frac{\vartheta'(\sqrt t)}{2\pi\sqrt t}\,dt + d\big[S(\sqrt t)\big] .
\end{equation}
Since $\vartheta'(T) = \tfrac12\Real\frac{\Gamma'}{\Gamma}(\tfrac14+\tfrac{iT}2) - \tfrac12\log\pi \approx \tfrac12\log(T/2\pi)$, the smooth part has density
\begin{equation}\label{eq:usmooth}
u_{\rm smooth}(t) = \frac{1}{4\pi\sqrt t}\Big(\Real\frac{\Gamma'}{\Gamma}\Big(\frac14+\frac{i\sqrt t}{2}\Big) - \log\pi\Big) \approx \frac{\log(\sqrt t/2\pi)}{4\pi\sqrt t} .
\end{equation}
This is $\tfrac12\log(\sqrt t/2\pi)$ times $1/(2\pi\sqrt t)$, the Thorin density of the positive stable law of index one half. So the Thorin measure of the reciprocal Xi function is a log corrected stable Thorin measure of index one half, plus the fluctuation $d[S(\sqrt t)]$.

The smooth part is not positive near the origin. The constant $1$ in the Riemann von Mangoldt formula is a unit of mass at $t = 0^+$, and $\vartheta'(\sqrt t) < 0$ for $\sqrt t < T_* = 6.2898\ldots$, so $u_{\rm smooth}$ is negative on $(0, 39.562)$ with total mass $\vartheta(T_*)/\pi = -1.1239\ldots$. The true Thorin measure has no mass at all below its first atom at $\gamma_1^2 = 199.79$, so on $(0, 199.79)$ the fluctuation must cancel both the unit mass and the negative density; indeed $S(\sqrt t) = -\vartheta(\sqrt t)/\pi - 1$ there.

\subsection*{The fluctuation and the primes}

Formally the fluctuation has density
\[
\frac{d}{dt}S(\sqrt t) = -\frac{1}{2\pi\sqrt t}\sum_{n\ge2}\frac{\Lambda(n)}{\sqrt n}\cos\big(\sqrt t\,\log n\big),
\]
a divergent series. It is the Thorin side of the convergent prime sum in \eqref{eq:weil}, since each cosine is the Laplace preimage of a Gaussian,
\[
\int_0^\infty e^{-xt}\cos(\sqrt t\,\log n)\frac{dt}{\sqrt{\pi t}} = \frac{e^{-(\log n)^2/(4x)}}{\sqrt x} .
\]
Putting the pieces together gives the signed Thorin density implied by the explicit formula,
\begin{equation}\label{eq:signedthorin}
U(dt) \;\text{``=''}\; \delta_{-1/4}(dt) + \frac{1}{4\pi\sqrt t}\Big(\Real\frac{\Gamma'}{\Gamma}\Big(\frac14+\frac{i\sqrt t}2\Big) - \log\pi\Big)dt - \frac{1}{2\pi\sqrt t}\sum_{n\ge2}\frac{\Lambda(n)}{\sqrt n}\cos(\sqrt t\log n)\,dt .
\end{equation}
The atom at $-1/4$ is the pole term $e^{x/4}$; it is where the unit of mass at $0^+$ really comes from. The equality is in quotation marks because the prime series is not even a tempered distribution in $t$: it holds only after Laplace transformation in $t$, where it becomes \eqref{eq:weil}.

\subsection*{Without RH}

By Theorem~\ref{thm:rh} there is then no positive Thorin measure at all, and the formulas above become formal. Two things survive, and they are the unconditional residue of the Thorin picture.

The first is the on line part. Let $U_0 = \sum_{\text{on line}}\delta_{\gamma^2}$. Then $\varphi = \varphi_0R$, where $\varphi_0(s) = \prod_{\text{on line}}(1+s/\gamma^2)^{-1}$ is a GGC Laplace transform without any hypothesis and $R(s) = \prod_{\text{off line pairs}}(1-s/s_\rho)^{-1}(1-s/\overline{s_\rho})^{-1}$ is positive on $[0,\infty)$; RH is the statement $R\equiv1$. Moreover the boundary argument counts exactly the zeros on the line, $\pi^{-1}\arg F(-t+i0) = N_0(\sqrt t)$, because an off line zero and its partner contribute conjugate factors whose arguments cancel on the cut.

The second is the tilted archimedean law. For $\lambda > 1/4$ let $\mu_\lambda(dy) = e^{-\lambda y}\mu(dy)/\varphi(\lambda)$ be the Esscher tilt of $\mu$; its Laplace transform is $\varphi(s+\lambda)/\varphi(\lambda)$ and its Lévy density $e^{-\lambda x}\Theta(x)/x$. Let $N$ have the zeta distribution of index $\sigma = \tfrac12+\sqrt\lambda$, and let $T_{\log N}$ be the first passage time of $\sqrt2B_t + 2\sqrt\lambda t$ to $\log N$; its Lévy density is $e^{-\lambda x}P(x)/x$, with $P$ the prime sum in \eqref{eq:weil}. Then $Y_\lambda + T_{\log N}$, with $Y_\lambda \sim \mu_\lambda$ independent, has the purely archimedean Laplace transform $\widetilde G(\tfrac12+\sqrt\lambda)/\widetilde G(\tfrac12+\sqrt{s+\lambda})$, where $\widetilde G(z) = \tfrac12z(z-1)\pi^{-z/2}\Gamma(z/2)$. Its canonical function has an exact signed Thorin representation,
\begin{equation}\label{eq:thorinR}
A(x)e^{-\lambda x} = e^{-(\lambda-1/4)x} + \int_\lambda^\infty e^{-zx}\,\frac{h(\sqrt{z-\lambda})}{4\pi\sqrt{z-\lambda}}\,dz, \qquad h(r) = \Real\frac{\Gamma'}{\Gamma}\Big(\frac14+\frac{ir}2\Big) - \log\pi,
\end{equation}
with a unit atom at $\lambda - \tfrac14$ and a density that is negative exactly on $(\lambda, \lambda+39.562)$. This is the smooth part \eqref{eq:usmooth} shifted by $\lambda$, with the atom of \eqref{eq:signedthorin} moved from $-1/4$ to $\lambda - 1/4 > 0$, and it is now an exact statement about a genuine infinitely divisible law.

\subsection*{A remark on normalisation}

The 2 sin$^2$ lemma of the next section writes the reciprocal as a function of $u = w^2/2$ instead of $s = w^2$. In that variable
\[
\frac{\xi(1/2)}{\xi(1/2+w)} = \exp\Big(-\int_0^\infty\big(1-e^{-w^2t/2}\big)\frac{\nu_{1/2}(t)}{t}\,dt\Big), \qquad \nu_{1/2}(t) = \Theta(t/2),
\]
and under RH the Thorin measure is $\sum_{\gamma>0}\delta_{\gamma^2/2}$. All statements of this section transfer by the scaling $t \mapsto t/2$.

\section{How the 2 sin\texorpdfstring{$^2$}{2} lemma turns a Lévy measure into a Thorin measure}\label{sec:sin2}

The two measures of this note live on different sides of a change of variable, and the device that carries one to the other is the 2 sin$^2$ lemma. This section explains how it works, in prose first and then in formulas.

\subsection*{The two sides and the statement}
On the primal side one has an even entire function $f$, such as $\xi$, and looks at $\log f(\alpha + w)$ as a function of $w$ at a base point $\alpha$. For $\xi$ at $\alpha > 1$ this function has a Lévy Khintchine representation of convex type,
\begin{equation}\label{eq:primal}
\log\frac{f(\alpha+w)}{f(\alpha)} - w\frac{f'(\alpha)}{f(\alpha)} = \int_0^\infty\big(e^{-wx} - 1 + wx\big)\,e^{-\alpha x}\,\frac{\mu(dx)}{x},
\end{equation}
with a measure $\mu$ assembled from the Euler product, the gamma factor and the pole: the Euler product gives $\mu = \sum_n\Lambda(n)\delta_{\log n}$, the gamma function gives a positive density, the pole factors $z(z-1)$ give negative mass. The right side of \eqref{eq:primal} is a positive combination, when $\mu \ge 0$, of the elementary convex functions $w \mapsto e^{-wx} - 1 + wx$. On the dual side one has the reciprocal $f(\alpha)/f(\alpha+w)$ written as a function of $u = w^2/2$, and asks whether it is the Laplace transform of a GGC, that is, whether its Laplace exponent in $u$ has the Thorin form $\int\log(1+u/z)\,U(dz)$ with $U \ge 0$. The lemma says that each elementary convex function on the primal side is a Thorin Bernstein function on the dual side, with an explicit nonnegative Thorin density, so that $\mu \ge 0$ transports to $U \ge 0$.

Here is the statement. For $x, w > 0$ and $u = w^2/2$,
\begin{equation}\label{eq:sin2}
e^{-wx} + wx - 1 = \int_0^\infty\log\Big(1 + \frac{u}{z}\Big)\,\tau_x(z)\,dz, \qquad \tau_x(z) = \frac{x}{\pi\sqrt{2z}}\;2\sin^2\!\Big(x\sqrt{\tfrac z2}\Big) \;\ge 0 .
\end{equation}
So the compensated exponential is, in the variable $u$, the Laplace exponent of a GGC whose Thorin density is a nonnegative oscillating function of $z$ with amplitude $x/(\pi\sqrt{2z})$. Integrating \eqref{eq:sin2} against $e^{-\alpha x}\mu(dx)/x$ turns \eqref{eq:primal} into
\[
\begin{gathered}
\log\frac{f(\alpha)}{f(\alpha+w)} + w\frac{f'(\alpha)}{f(\alpha)} = -\int_0^\infty\log\Big(1+\frac uz\Big)\tau_\alpha(z)\,dz,\\
\tau_\alpha(z) = \frac{1}{\pi\sqrt{2z}}\int_0^\infty 2\sin^2\!\Big(x\sqrt{\tfrac z2}\Big)e^{-\alpha x}\mu(dx),
\end{gathered}
\]
which is the Thorin representation of the reciprocal, with the linear term $w f'(\alpha)/f(\alpha)$ removed by the factor $e^{w f'/f}$. If $\mu \ge 0$ then $\tau_\alpha \ge 0$ and the reciprocal is a GGC Laplace transform in $u$. For $\xi$ at $\alpha = 1/2$ the linear term vanishes because $\xi'(1/2) = 0$, and the Thorin measure, when it exists, is $\sum\delta_{\gamma^2/2}$.

\subsection*{Why the lemma is true: two steps}
The proof has two steps, and each has a probabilistic meaning.

Step one is subordination. The exponential $e^{-wx}$ is the Laplace transform, in the variable $u = w^2/2$, of the first passage time $T_x$ of a standard Brownian motion to the level $x$, whose density is $x(2\pi t^3)^{-1/2}e^{-x^2/(2t)}$. The linear term $wx$ is the Laplace exponent, in $u$, of the positive stable subordinator of index one half, whose Lévy density is $x(2\pi t^3)^{-1/2}$; indeed $\int_0^\infty(1-e^{-ut})\,x(2\pi t^3)^{-1/2}dt = x\sqrt{2u} = wx$. Subtracting,
\begin{equation}\label{eq:stepA}
e^{-wx} + wx - 1 = \int_0^\infty\big(1 - e^{-ut}\big)\,k_x(t)\,dt, \qquad k_x(t) = \frac{x}{\sqrt{2\pi t^3}}\Big(1 - e^{-x^2/(2t)}\Big) \ge 0 .
\end{equation}
So the compensated exponential is a Bernstein function of $u$, with Lévy density $k_x$: the stable Lévy density thinned by the factor $1 - e^{-x^2/(2t)}$, which is the probability that a Brownian bridge of duration $t$ stays below $x/2$. The convexity of $e^{-wx}-1+wx$ in $w$ has become positivity of a Lévy density in $u$; this is the step that explains why the $\sqrt s$ variable is the right one for the reciprocal.

Step two is the sine squared identity. A Bernstein function with Lévy density $k(t)$ is a Thorin Bernstein function exactly when $t\,k(t)$ is completely monotone, that is, a Laplace transform $t\,k(t) = \int e^{-tz}\tau(z)\,dz$ with $\tau \ge 0$, and then $\tau$ is the Thorin density, because
\[
\int_0^\infty(1 - e^{-ut})\,k(t)\,dt = \int_0^\infty(1-e^{-ut})\int_0^\infty e^{-tz}\tau(z)\,dz\,\frac{dt}{t} = \int_0^\infty\log\Big(1+\frac uz\Big)\tau(z)\,dz
\]
by Frullani's integral. For the thinned stable density $t\,k_x(t) = x(2\pi)^{-1/2}\,t^{-1/2}(1 - e^{-x^2/(2t)})$, and the two pieces are both Laplace transforms:
\[
\frac{1}{\sqrt t} = \int_0^\infty e^{-tz}\frac{dz}{\sqrt{\pi z}}, \qquad \frac{e^{-x^2/(2t)}}{\sqrt t} = \int_0^\infty e^{-tz}\cos\big(x\sqrt{2z}\big)\frac{dz}{\sqrt{\pi z}},
\]
the first by the gamma integral and the second by the Fourier transform of a Gaussian after the substitution $z = r^2/2$. Subtracting, $1 - \cos(2\theta) = 2\sin^2\theta$ gives
\begin{equation}\label{eq:stepB}
\frac{1 - e^{-x^2/(2t)}}{\sqrt t} = \int_0^\infty e^{-tz}\,2\sin^2\!\Big(x\sqrt{\tfrac z2}\Big)\frac{dz}{\sqrt{\pi z}},
\end{equation}
and \eqref{eq:sin2} follows with $\tau_x(z) = x(2\pi)^{-1/2}\cdot 2\sin^2(x\sqrt{z/2})(\pi z)^{-1/2} = x\,2\sin^2(x\sqrt{z/2})/(\pi\sqrt{2z})$. The point of the identity is the sign: the Laplace density of the stable piece is $1/\sqrt{\pi z}$, the Laplace density of the first passage piece is $\cos(x\sqrt{2z})/\sqrt{\pi z}$, and their difference is nonnegative because $1 - \cos \ge 0$. A Lévy density of the form stable minus first passage is therefore always Thorin, and this is what the lemma exploits. Both \eqref{eq:stepA} and \eqref{eq:stepB} were checked numerically to ten digits, and \eqref{eq:sin2} itself to four digits with an oscillatory quadrature.

\subsection*{A worked example: the gamma function}
Take $f(w) = \Gamma(w)$, for which Gauss's formula gives \eqref{eq:primal} with the positive density $e^{-\alpha x}\mu(dx)/x = e^{-\alpha x}dx/(x(1-e^{-x}))$, so that $\mu(dx) = dx/(1-e^{-x}) = \sum_{k\ge0}e^{-kx}dx$. The sine squared transform can be done term by term,
\[
\int_0^\infty 2\sin^2\!\Big(x\sqrt{\tfrac z2}\Big)e^{-(\alpha+k)x}dx = \frac{1}{\alpha+k} - \frac{\alpha+k}{(\alpha+k)^2+2z} = \frac{2z}{(\alpha+k)\big((\alpha+k)^2+2z\big)},
\]
so the Thorin density of $\Gamma(\alpha)/\Gamma(\alpha+\sqrt{2u})\,e^{\sqrt{2u}\,\Gamma'(\alpha)/\Gamma(\alpha)}$ is
\[
\tau^\Gamma_\alpha(z) = \frac{\sqrt{2z}}{\pi}\sum_{k\ge0}\frac{1}{(\alpha+k)\big((\alpha+k)^2+2z\big)},
\]
a positive density behaving like $\sqrt z$ at the origin and like $\pi^{-1}(2z)^{-1/2}\,\Gamma'(\alpha)/\Gamma(\alpha)$ type terms at infinity. At $\alpha = 1$, $u = 0.7$ the Thorin integral $\int\log(1+u/z)\tau^\Gamma_1(z)dz = 0.7708$ agrees with $\log\Gamma(1+\sqrt{1.4}) - \sqrt{1.4}\,\Gamma'(1)/\Gamma(1) = 0.7709$. This is the archimedean part of the Thorin density of the reciprocal Xi function in Section 4, in the $u$ normalisation and with the pole factors not yet included.

\subsection*{Where it fails for $\xi$, and what replaces it}
For $\xi(z) = \tfrac12z(z-1)\pi^{-z/2}\Gamma(z/2)\zeta(z)$ at a base point $\alpha > 1$ the measure $\mu$ in \eqref{eq:primal} is a sum of four parts: the Euler product, positive atoms at $\log n$; the gamma factor, a positive density; the power of $\pi$, a pure linear term absorbed in $f'/f$; and the two pole factors $z$ and $z-1$, whose logarithms are concave in $w$ and so contribute negative mass. The lemma transports each part separately, and the negative parts produce a negative Thorin density concentrated near $z = 0$, of size $-\sqrt{2z}/(\pi c(c^2+2z))$ for a pole factor $w + c$. So positivity of the Thorin density of the reciprocal $\xi$ at $\alpha>1$ is not a termwise consequence of the lemma but a cancellation, which is nonetheless established in \cite{Ladder}. At $\alpha = 1/2$ the situation is worse: the Euler product no longer converges, the prime part of $\tau_\alpha$ becomes the divergent constant $(\pi\sqrt{2z})^{-1}\sum\Lambda(n)n^{-1/2}$ minus a cosine series that is not a tempered distribution, and the cancellation can no longer be organised term by term. This is exactly the point at which the Lévy side of this note takes over: the Laplace transform in $t$ of the divergent Thorin side is the convergent explicit formula \eqref{eq:weil}, the Thorin constant becomes the pole term $e^{x/4}$ and the archimedean terms, and the cosine series becomes the Gaussian prime sum. In the language of the lemma, at $\alpha = 1/2$ step one still works, giving the Lévy density $\Theta(x)/x$ with $\Theta$ positive and nonincreasing, while step two, the complete monotonicity of $\Theta$ that would produce a positive Thorin density, is the Riemann Hypothesis.

\section{What is unconditional and what is RH}\label{sec:division}

This section separates the two kinds of statement precisely and gives the mechanism behind each, with proofs or proof sketches. The dictionary of Section~\ref{sec:law} reduces everything to the sign of $T_k(x) = (-1)^k\Theta^{(k)}(x)$ on a range of $(k,x)$, so the question is always: for which $(k,x)$ can a zero off the line make a term negative, and does anything cancel it? Throughout, $H = 3\cdot10^{12}$ is the height to which all zeros are known to lie on the line~\cite{PlattTrudgian}, $\rho = \tfrac12+\delta+i\gamma$ runs over zeros in the upper half plane, $|\delta| < \tfrac12$, and $T_k(x) = (-1)^k\Theta^{(k)}(x) = \sum_{\gamma>0}\gamma_\rho^{2k}e^{-\gamma_\rho^2x}$ with $\gamma_\rho = \gamma - i\delta$.

\subsection*{Positivity: the phase and the window}
Each term of $T_k(x)$ has modulus $|\gamma_\rho|^{2k}e^{-(\gamma^2-\delta^2)x}$ and phase
\begin{equation}\label{eq:phase}
\Psi_\rho(k,x) = 2x\gamma\delta - 2k\arctan(\delta/\gamma) = 2x\delta\,\frac{\gamma^2 - T_0^2}{\gamma} + O\Big(\frac{k\delta^3}{\gamma^3}\Big), \qquad T_0 = \sqrt{k/x}.
\end{equation}
A term can have negative real part only if $|\Psi_\rho| > \pi/2$. Every off line zero has $\gamma > H$, so for $k \le 10^{10}$ the term $2k\arctan(\delta/\gamma)$ in \eqref{eq:phase} is below $0.004$, while $|2x\gamma\delta| < x\gamma$ because $|\delta| < \tfrac12$; hence an off line term is nonnegative unless $\gamma > 0.78/x$, a deliberately conservative threshold. Such zeros carry weight at most $e^{-\gamma^2x} \le e^{-0.78\gamma}$, which the explicit zero counting bound shows is negligible against the first zero's term $\gamma_1^{2k}e^{-\gamma_1^2x}$. This gives $T_k > 0$ for all $x$ and all $k \le 10^{10}$, hence (a) and (b): the law is infinitely divisible and self decomposable, unconditionally.

Reading \eqref{eq:phase} the other way, with $T_0 = \sqrt{k/x}$ the centre of the Gaussian window $f(T) = T^{2k}e^{-xT^2}$ and $1/(2\sqrt x)$ its width, the phase of a zero at distance $d$ from the centre is about $4x\delta d$, so across one window width it changes by $2\delta\sqrt x$. For $x \le 1/20$ and $|\delta| < 1/2$ this is below $0.22$, every zero within the window has positive real part, zeros far outside are Gaussian suppressed, and a window of length $4$ above $T_0$ contains at least $0.12\log T_0$ zeros by the explicit counting bound. This gives $T_k > 0$ for all $k$ when $x \le 1/20$. Combining with the first argument, a failure needs $k > 10^{10}$, $x > 1/20$ and $T_0 > H/2$, hence $k = xT_0^2 > H^2/80$. So the ladder holds to order $H^2/80$ at every $x$, and to every order at coarse resolution.

\subsection*{Positivity on the cut: the counting mechanism}
With $F(s) = \prod_\rho(1 - s/s_\rho)$, the argument of $F$ along the upper edge of the cut is a sum over zeros. An on line zero contributes $\arg(1 + (-t+i0)/\gamma^2)$, which jumps from $0$ to $\pi$ as $t$ passes $\gamma^2$. An off line zero and its partner give the two factors $(1-s/s_\rho)(1-s/\overline{s_\rho}) = |1 - s/s_\rho|^2 > 0$ for real $s$, with argument identically $0$. Hence
\begin{equation}\label{eq:count}
\frac1\pi\arg F(-t+i0) = N_0(\sqrt t) \quad\text{unconditionally}, \qquad = N(\sqrt t) \quad\text{under RH},
\end{equation}
where $N_0$ counts zeros on the line and $N$ all zeros. The same product gives the factorisation $\varphi = \varphi_0 R$: $\varphi_0 = \prod_{\rm on}(1+s/\gamma^2)^{-1}$ is a GGC Laplace transform for every zero set, and $R = \prod_{\rm off}|1-s/s_\rho|^{-2}$ is positive on $[0,\infty)$, has no pole on the cut, and is identically $1$ exactly under RH.

\subsection*{Identities that need no hypothesis}
\emph{The explicit formula.} The Weil explicit formula with the test function $e^{-r^2x}$ is an identity between the zero side $\sum_\rho e^{(\rho-1/2)^2x}$, which is $2\Theta(x)$, and the prime side \eqref{eq:weil}, valid for every $x > 0$ with no hypothesis, because the Gaussian is entire and decays in every strip. So \eqref{eq:weil} and its level form \eqref{eq:levelform}, the Lévy density \eqref{eq:levydensity}, the small time asymptotic \eqref{eq:small}, and the cumulant formula $\kappa_n = (n-1)!\sum\gamma_\rho^{-2n}$ are all unconditional; under RH the last sum is over real $\gamma$ and all cumulants are positive, while unconditionally the complex terms pair into real values, and every cumulant is positive, because $\kappa_n = (n-1)!\int_0^\infty x^{n-1}\Theta(x)\,dx$ and $\Theta > 0$.

\emph{Operators.} The multiplier $\varphi(|\xi|^2)^\theta$ on $\R^d$ is a Schwartz function for $\theta > 0$, because $F(\lambda) = \sum c_n\lambda^n$ with $c_n > 0$ by Riemann's representation of $\xi(1/2+w)$ as $8\int\Phi(u)\cosh(2wu)\,du$ with $\Phi > 0$. By Bochner the kernel of $\varphi(-\Delta)^\theta$ is nonnegative iff $\varphi^\theta(|\xi|^2)$ is positive definite, and by Schoenberg this holds in every dimension iff $\varphi^\theta$ is completely monotone for all $\theta>0$, iff $\varphi$ is infinitely divisible, iff $\Theta \ge 0$. So the probability kernels are unconditional, because (a) is.

\emph{The Esscher tilt.} For $\lambda > 1/4$ the Esscher tilt $\mu_\lambda$ has Lévy density $e^{-\lambda x}\Theta(x)/x$, and the explicit formula splits $\Theta = A - P$ with $P$ the prime sum. The identity
\[
\frac{\zeta(\tfrac12+\sqrt{s+\lambda})}{\zeta(\tfrac12+\sqrt\lambda)} = \E\,e^{-sT_{\log N}}, \qquad N \sim \text{Zeta}(\tfrac12+\sqrt\lambda),
\]
is the Euler product read as a first passage time: $\E\,e^{-sT_a} = e^{-a(\sqrt{s+\lambda}-\sqrt\lambda)}$ for the drifted Brownian motion, and $\E\,N^{-v} = \zeta(\sigma+v)/\zeta(\sigma)$. It needs only $\sigma > 1$. Multiplying Laplace transforms, $Y_\lambda + T_{\log N}$ has the archimedean law with exponent $\log\widetilde G(\tfrac12+\sqrt{s+\lambda}) - \log\widetilde G(\tfrac12+\sqrt\lambda)$, whose Lévy density $A(x)e^{-\lambda x}/x$ is positive because $A = \Theta + P$ with both terms positive, and whose Thorin representation \eqref{eq:thorinR} is read off from $A = e^{x/4} + \frac1{4\pi}\int e^{-r^2x}h(r)\,dr$ by substituting $z = r^2 + \lambda$. All of this is unconditional. What RH adds is that $\mu_\lambda$ is a GGC with Thorin measure $\sum\delta_{\gamma^2+\lambda}$, so that the arithmetic factor $T_{\log N}$ exactly cancels the negative Thorin mass $-1.124$ of the archimedean law on $(\lambda, \lambda+39.56)$ and leaves atoms only at $\gamma^2 + \lambda \ge 199.79 + \lambda$.

\subsection*{What requires RH}
The Thorin measure is defined by Stieltjes inversion \eqref{eq:inversion}, and \eqref{eq:count} shows that its distribution function along the cut is $N_0(\sqrt t)$ unconditionally. But a positive measure with this distribution function has Laplace transform $\Theta_0(x) = \sum_{\rm on}e^{-\gamma^2x}$, not $\Theta$, and the two differ by $\Theta_{\rm off}(x) = \sum_{\rm off}e^{(\delta^2-\gamma^2)x}\cos(2\delta\gamma x)$, which is the Laplace transform of no positive measure: its atoms sit at the complex points $-s_\rho$. So the statement ``$\Theta$ is the Laplace transform of the measure with distribution function $N(\sqrt t)$'' is exactly RH, and the formulas \eqref{eq:thorinRH} to \eqref{eq:usmooth} are RH statements. The signed formula \eqref{eq:signedthorin} is unconditional only as a formal Laplace preimage of \eqref{eq:weil}.

\emph{The explicit representations} are RH statements for the same reason. $\mu = $ law of $\sum E_\gamma/\gamma^2$ requires real $\gamma$. The driver $D$ as a superposition of compound Poisson processes with $\mathrm{Exp}(\gamma^2)$ jumps requires $-\Theta' = \sum\gamma^2e^{-\gamma^2x}$ with real $\gamma$; unconditionally $-\Theta'$ is positive but its terms are complex. The Krein string with spectrum $\{\gamma^2\}$ requires RH and simple zeros, by Kent's characterisation of hitting times~\cite{Kent}. The bound $S(T) = O(\log T/\log\log T)$, used for the Kolmogorov distance between $U$ and its smooth part, is Littlewood's conditional estimate~\cite{Titchmarsh}; unconditionally $S(T) = O(\log T)$.

\subsection*{The converse direction and the order of monotonicity}
If RH fails, the ladder breaks at finite order. With $w_* = -s_\rho$ for an off line zero maximising $\sec\arg w$, the Post Widder transform $\frac{(-1)^k}{k!}(k/y)^{k+1}\Theta^{(k)}(k/y)$ at $y = \Real w_*$ grows like $(\sec\arg w_*)^k\cos(k\eta)$ with $\eta = \tan|\arg w_*| - |\arg w_*| > 0$, so it is negative for infinitely many $k$. For an isolated zero the Gaussian window argument gives a failure already at order $O(\gamma^2\log\gamma/\delta)$. So $K_* = \sup\{K : T_k \ge 0 \text{ for all } k \le K, x > 0\}$ satisfies $H^2/80 \le K_*$ unconditionally, and $K_* < \infty$ exactly when RH fails. The order of monotonicity of the canonical function is, up to logarithms, the square of the height of the first off line zero.

\subsection*{The division in one table}
Table~\ref{tab:division} collects the statements of this section with the mechanism behind each and its status.

\begin{table}[tbp]
\centering
\small
\begin{tabular}{p{4.6cm}p{4.0cm}p{5.0cm}}
\toprule
Statement & Mechanism & Status\\
\midrule
$\Theta > 0$, $-\Theta' > 0$; infinitely divisible, self decomposable & phase \eqref{eq:phase}, zeros below $H$, counting bound & unconditional\\
$T_k > 0$ for $k \le H^2/80$, and for all $k$ at $x \le 1/20$ & Gaussian window, phase across one width $2\delta\sqrt x$ & unconditional\\
Explicit formula forms of $\Theta$, Lévy density, cumulants, asymptotics & Weil formula with Gaussian test function & unconditional\\
Boundary argument $= N_0(\sqrt t)$; $\varphi = \varphi_0R$ & conjugate pairs cancel on the cut & unconditional\\
Probability kernels of $\varphi(-\Delta)^\theta$ & Bochner, Schoenberg, $c_n > 0$ & unconditional\\
Tilted balance and archimedean Thorin measure & Euler product as first passage; $A = \Theta + P > 0$ & unconditional ($\lambda > 1/4$)\\
Ladder breaks at finite order if RH fails & Post Widder growth $(\sec\arg w)^k$ & unconditional converse\\
\midrule
$\Theta$ completely monotone; $\mu$ a GGC & poles of $\psi'$ on the cut & equivalent to RH\\
Thorin measure $\sum\delta_{\gamma^2}$, distribution $N(\sqrt t)$ & Stieltjes inversion & RH\\
$\mu = \sum E_\gamma/\gamma^2$; exponential jump driver; Krein string & real ordinates & RH (string also needs simple zeros)\\
$S(T) = O(\log T/\log\log T)$ & Littlewood & RH\\
\bottomrule
\end{tabular}
\caption{What is unconditional and what is equivalent to RH.}
\label{tab:division}
\end{table}

\section{Weak convergence and Bondesson's continuity theorem}\label{sec:weak}

The GGC property is a closed condition, and this is the precise sense in which a limiting argument could finish a proof of RH. We state the continuity theorem in the form needed, prove it from the dictionary of Section~\ref{sec:law}, and then explain why closure alone cannot supply the missing input.

\begin{theorem}[Continuity of GGC laws]\label{thm:continuity}
Let $\mu_n$ be GGC laws on $[0,\infty)$ with Laplace exponents $\psi_n(s) = a_ns + \int\log(1+s/t)\,U_n(dt)$, and suppose $\mu_n \to \mu$ weakly, where $\mu$ is a probability law on $[0,\infty)$. Then $\mu$ is a GGC. Its exponent $\psi$ satisfies $\psi_n \to \psi$ and $\psi_n' \to \psi'$ pointwise on $(0,\infty)$, and its Thorin measure $U$ is the vague limit on $(0,\infty)$ of the Thorin measures $U_n$, that is $\int g\,dU_n \to \int g\,dU$ for every continuous $g$ with compact support in $(0,\infty)$; mass of $U_n$ escaping to infinity contributes to the drift of the limit. Conversely, if $\psi_n$ converges pointwise on $(0,\infty)$ to a function $\psi$ with $\psi(0+) = 0$, then $\mu_n$ converges weakly to the GGC law with exponent $\psi$.
\end{theorem}

\begin{proof}
By the continuity theorem for Laplace transforms, weak convergence on $[0,\infty)$ is equivalent to pointwise convergence $\varphi_n \to \varphi$ on $(0,\infty)$, hence $\psi_n = -\log\varphi_n \to \psi$ pointwise. Each $\psi_n$ is a Bernstein function, in particular concave on $(0,\infty)$, and pointwise convergence of concave functions implies convergence of their derivatives at every point where the limit is differentiable; $\psi$ is smooth on $(0,\infty)$, being minus the logarithm of a Laplace transform, so $\psi_n' \to \psi'$ everywhere on $(0,\infty)$. Each $\psi_n'$ is a Stieltjes function $a_n + \int U_n(dt)/(t+s)$, and the cone of Stieltjes functions is closed under pointwise convergence on $(0,\infty)$ \cite[Chapter 7]{SSV}; hence $\psi' = a + \int U(dt)/(t+s)$ with $a \ge 0$ and $U \ge 0$, and integrating, $\psi(s) = as + \int\log(1+s/t)U(dt)$, so $\mu$ is a GGC. The representing measures of a pointwise convergent sequence of Stieltjes functions converge vaguely on $(0,\infty)$, by the Stieltjes inversion formula together with the uniform bound $\int(1+t)^{-1}U_n(dt) \le \psi_n'(1) \to \psi'(1)$, which gives tightness on compact subsets of $(0,\infty)$; mass may escape to $t = \infty$, where $(t+s)^{-1} \approx t^{-1}$ no longer depends on $s$, so that it appears in the limit as a constant term of $\psi'$, that is, as drift. For the converse, $e^{-\psi_n} \to e^{-\psi}$ pointwise with $e^{-\psi(0+)} = 1$, so by Laplace continuity $\mu_n$ converges weakly to a probability law with exponent $\psi$, and by the first part that law is a GGC.
\end{proof}

The theorem is Bondesson's \cite{Bondesson}, and the proof shows why it holds: weak convergence of the laws is pointwise convergence of the Stieltjes functions $\psi_n'$, and the Thorin measure is the boundary value of $\psi'$ on the cut, so it moves continuously with the laws. Applied to the reciprocal Xi law it gives a clean criterion.

\begin{corollary}\label{cor:closed}
If $\mu$ is the weak limit of a sequence of GGC laws, then RH holds. Equivalently, the set of RH statements expressible through the reciprocal Xi law is closed under weak limits.
\end{corollary}

So a proof of RH could take the form: construct GGC laws $\mu_n$ converging weakly to $\mu$. The difficulty is entirely in the construction. Every natural family converging to $\mu$ is a family of GGCs near $\mu$ only under a hypothesis equivalent to RH, as Table~\ref{tab:families} shows.

\begin{table}[htbp]
\centering
\small
\begin{tabular}{p{3.0cm}p{4.6cm}p{6.0cm}}
\toprule
Family & GGC when & What closure gives\\
\midrule
Heat flow $\mu_t$, from $H_t$ & $t \ge \Lambda$ (de Bruijn Newman~\cite{deBruijn,Newman}) & $\mu_\Lambda$ is a GGC; reaching $\mu = \mu_0$ needs $\Lambda \le 0$, which with $\Lambda \ge 0$ is RH\\
Base point $\mu^{(\alpha)}$, Section~\ref{sec:base} & $\alpha \ge \sigma_{\max} = \sup\Real\rho$ & $\mu^{(\alpha)}$ is a GGC for $\alpha \ge 1$ unconditionally; reaching $\alpha = 1/2$ needs $\sigma_{\max} = 1/2$, which is RH\\
On line part $\varphi_0$ & always & no limit is involved; $\varphi_0 = \varphi$ is RH\\
Truncated Euler products & not applicable & they do not converge on the critical line, so give no weak limit at $\alpha = 1/2$\\
Random matrix models & always (zeros on a circle) & they match $\zeta$ in zero statistics, not as functions, so say nothing about $\mu$\\
\bottomrule
\end{tabular}
\caption{Approximating families for the reciprocal Xi law.}
\label{tab:families}
\end{table}

The reason is structural. A closed condition can be reached by a limit only from approximants that already satisfy it arbitrarily close to the target, and for every family in the table that is a restatement of RH. The base point family is the most promising of the three genuine ones, because its GGC property is unconditional on a half line and is proved by arithmetic, so the next section develops it in full.

\section{The laws at base point \texorpdfstring{$\alpha$}{alpha} and their descent into the critical strip}\label{sec:base}

\subsection*{The product formula and the base point theorem}

For real $\alpha$ write each zero in the upper half plane as $\rho = \beta + i\gamma$, and put
\[
c_\rho = \alpha - \beta, \qquad q_\rho = c_\rho^2 + \gamma^2, \qquad \sigma_{\max} = \sup_\rho\beta .
\]
Unconditionally $\tfrac12 \le \sigma_{\max} \le 1$, with no zero on $\Real\rho = 1$, but $\sigma_{\max} = 1$ is not excluded by any known zero free region, since these shrink toward the line $\Real\rho = 1$ as $\gamma\to\infty$. RH is $\sigma_{\max} = \tfrac12$.

\begin{theorem}[Product formula]\label{thm:product}
For every real $\alpha$ and every $s \ge 0$ (and, by analytic continuation, wherever both sides are finite),
\[
\frac{\xi(\alpha)}{\xi(\alpha+\sqrt s)} = \prod_{\gamma>0}\frac{q_\rho}{(c_\rho+\sqrt s)^2 + \gamma^2},
\]
the product running over zeros in the upper half plane with multiplicity, and converging absolutely.
\end{theorem}

\begin{proof}
Hadamard's product $\xi(z) = \xi(0)\prod_\rho(1 - z/\rho)$, with $\rho$ and $1-\rho$ paired, gives $\xi(\alpha+w)/\xi(\alpha) = \prod_\rho(1 + w/(\alpha-\rho))$. Grouping $\rho$ with $\overline\rho$, $(1+w/(\alpha-\rho))(1+w/(\alpha-\overline\rho)) = ((c_\rho+w)^2+\gamma^2)/q_\rho = 1 + (2c_\rho w + w^2)/q_\rho$, and $\sum q_\rho^{-1} < \infty$, so the grouped product converges absolutely and no exponential factors are needed. Put $w = \sqrt s$.
\end{proof}

At $\alpha = \tfrac12$ under RH this is $F(s)^{-1} = \prod(1+s/\gamma^2)^{-1}$. Numerically, with the first $600$ zeros and a Riemann von Mangoldt tail, the product agrees with $\xi(\alpha)/\xi(\alpha+w)$ to six digits at $\alpha = 1$ and $0.75$ and $w = 0.5, 1, 3$; for example $0.954929659$ against $0.954929048$ at $\alpha = 1$, $w = 1$.

\begin{theorem}[Base point]\label{thm:base}
The function $\varphi^{(\alpha)}(s) = \xi(\alpha)/\xi(\alpha+\sqrt s)$ is the Laplace transform of a GGC law $\mu^{(\alpha)}$ if and only if $\alpha \ge \sigma_{\max}$. In that case the drift is zero and the Thorin measure is $U_\alpha = \sum_{\gamma>0}u_{\rho,\alpha}(z)\,dz$, with atoms $\delta_{\gamma^2}$ for zeros with $c_\rho = 0$, where
\begin{equation}\label{eq:lorentz}
u_{\rho,\alpha}(z) = \frac{c_\rho\,(z + q_\rho)}{\pi\sqrt z\,\big[(z - q_\rho)^2 + 4c_\rho^2z\big]} \ge 0 .
\end{equation}
Its distribution function is a sum of arctangents,
\begin{equation}\label{eq:arctan}
U_\alpha\big((0,T^2]\big) = \frac1\pi\sum_{\gamma>0}\Big[\arctan\frac{T-\gamma}{c_\rho} + \arctan\frac{T+\gamma}{c_\rho}\Big] = \frac1\pi\Delta_{0}^{T}\arg\xi(\alpha+i\tau),
\end{equation}
the last term being the continuous change in the argument of $\xi$ along the vertical line $\Real z = \alpha$ from height $0$ to height $T$. Unconditionally, $\mu^{(\alpha)}$ is a GGC for every $\alpha \ge 1$.
\end{theorem}

\begin{proof}
Each factor $q/((c+\sqrt s)^2+\gamma^2) = [(1+\sqrt s/a)(1+\sqrt s/\overline a)]^{-1}$ with $a = c + i\gamma$ has exponent $\psi_\rho(s) = \log((c+\sqrt s)^2+\gamma^2) - \log q$, and
\[
\psi_\rho'(s) = \frac{c+\sqrt s}{\sqrt s\,\big((c+\sqrt s)^2+\gamma^2\big)} = \frac12\Big[\frac{1}{\sqrt s(\sqrt s+a)} + \frac{1}{\sqrt s(\sqrt s+\overline a)}\Big].
\]
Stieltjes inversion with $\sqrt{-z+i0} = i\sqrt z$ gives the density
\[
-\frac1\pi\Imag\psi_\rho'(-z+i0) = \frac1\pi\Real\frac{a}{\sqrt z\,(z+a^2)},
\]
and with $a^2 = c^2-\gamma^2+2ic\gamma$ the real part simplifies to \eqref{eq:lorentz}, using $(z+c^2-\gamma^2)^2 + 4c^2\gamma^2 = (z-q)^2 + 4c^2z$. It is nonnegative exactly when $c \ge 0$, and $\psi_\rho'(s) \to 0$ as $s \to \infty$, so the drift vanishes. Summing over zeros, which converges by Theorem~\ref{thm:product}, gives the Thorin representation when every $c_\rho \ge 0$. Integrating \eqref{eq:lorentz} in $z = \tau^2$ gives the arctangents, since $\frac{d}{d\tau}\big[\arctan\frac{\tau-\gamma}{c} + \arctan\frac{\tau+\gamma}{c}\big] = \frac{2c(\tau^2+q)}{(\tau^2-q)^2+4c^2\tau^2}$, and each arctangent pair is the change of $\arg(\alpha + i\tau - \rho) + \arg(\alpha+i\tau-\overline\rho)$, whose sum over zeros is the change of $\arg\xi(\alpha+i\tau)$ by the product formula. If some $c_\rho < 0$, the factor of that zero has poles at $\sqrt s = -c_\rho \pm i\gamma$, which have positive real part, so at $s = (|c_\rho| \pm i\gamma)^2$, off the negative axis; a GGC Laplace transform is holomorphic on $\C\setminus(-\infty,0]$, so $\varphi^{(\alpha)}$ is not one. Finally every zero has $\beta < 1$, since $\zeta$ has no zeros with $\Real z \ge 1$, so $c_\rho > 0$ for all $\rho$ when $\alpha \ge 1$.
\end{proof}

The Lévy measure of $\mu^{(\alpha)}$ is equally explicit when $\alpha \ge \sigma_{\max}$. Since $\int_0^\infty e^{-sx}e^{a^2x}\operatorname{erfc}(a\sqrt x)\,dx = 1/(\sqrt s(\sqrt s+a))$ for $\Real a > 0$, the canonical function is
\begin{equation}\label{eq:erfc}
k_\alpha(x) = \sum_{\gamma>0}\Real\Big[e^{(c_\rho+i\gamma)^2x}\operatorname{erfc}\big((c_\rho+i\gamma)\sqrt x\big)\Big], \qquad \nu_\alpha(dx) = \frac{k_\alpha(x)}{x}\,dx,
\end{equation}
checked against $\psi_\rho'$ by quadrature for the first zero at $c = 0.5$, $s = 3$ to $25$ digits. At $\alpha = 1/2$ under RH, $a = i\gamma$ and $\Real[e^{-\gamma^2x}\operatorname{erfc}(i\gamma\sqrt x)] = e^{-\gamma^2x}$, since $\operatorname{erf}(iy)$ is imaginary, recovering $\Theta$. For $\alpha > 1/2$ the law is heavy tailed. The exponent behaves like $\psi(s) \sim \mathsf B_\alpha\sqrt s$ as $s \to 0$, with
\[
\mathsf B_\alpha = \sum_{\gamma>0}\frac{2c_\rho}{q_\rho} = \frac{\xi'(\alpha)}{\xi(\alpha)},
\]
so by the Tauberian theorem $\mu^{(\alpha)}((y,\infty)) \sim \mathsf B_\alpha/\sqrt{\pi y}$ and $\mu^{(\alpha)}$ has infinite mean. Numerically $\mathsf B_1 = 0.023096$, $\mathsf B_{0.75} = 0.011551$, $\mathsf B_{0.6} = 0.004621$, $\mathsf B_{0.52} = 0.000924$, and $\mathsf B_{1/2} = 0$ because $\xi'(1/2) = 0$, which is why $\mu = \mu^{(1/2)}$ alone has a finite mean.

\begin{remark}[The compensated normalisation]
The 2 sin$^2$ lemma of Section~\ref{sec:sin2} uses the compensated reciprocal $\varphi^{(\alpha)}(s)\,e^{\mathsf B_\alpha\sqrt s}$, which has finite mean but whose Thorin density is $U_\alpha - \mathsf B_\alpha/(2\pi\sqrt z)$, since $e^{-\mathsf B\sqrt s}$ is the Laplace transform of a positive stable law with Thorin density $\mathsf B/(2\pi\sqrt z)$. Positivity of the compensated density is therefore a cancellation, as Section~\ref{sec:sin2} found from the prime side. The plain reciprocal of Theorem~\ref{thm:base} has no such subtraction: from the zero side its Thorin density is a sum of nonnegative Lorentzians as soon as all zeros lie to the left of $\alpha$.
\end{remark}

\subsection*{Descent as a backward Cauchy flow}

\begin{figure}[htbp]
\centering
\includegraphics[width=0.9\textwidth]{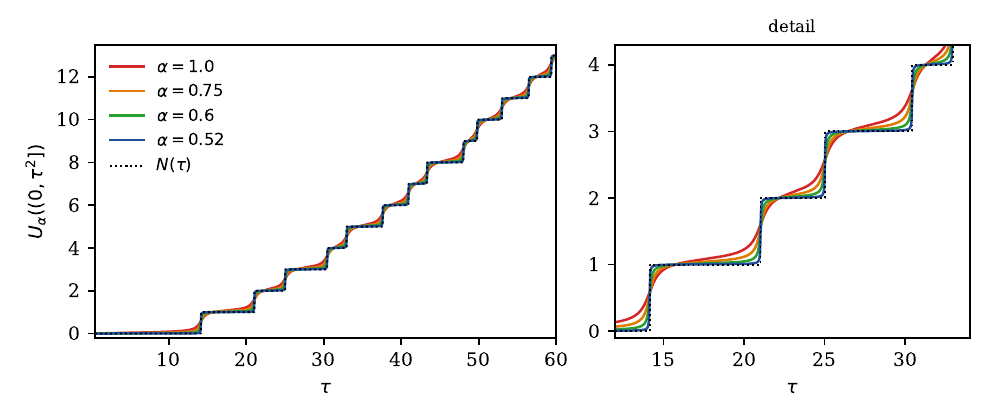}
\caption{The backward Cauchy flow. The Thorin distribution functions $U_\alpha((0,\tau^2])$ of \eqref{eq:arctan} at $\alpha = 1$ (red), $0.75$ (orange), $0.6$ (green) and $0.52$ (blue) sharpen onto the zero counting function $N(\tau)$ (dotted) as the base point descends; the right panel shows the first four zeros. By Theorem~\ref{thm:continuity} the corresponding laws converge weakly to the reciprocal Xi law, and every curve is nondecreasing because every zero shown lies on the line.}
\label{fig:descent}
\end{figure}

Formula \eqref{eq:arctan} says that, in the ordinate variable $\tau = \sqrt z$, the Thorin measure at base point $\alpha$ is the zero counting measure smoothed zero by zero with a Cauchy kernel of width $c_\rho = \alpha - \beta_\rho$. For zeros on the line all widths equal $\alpha - \tfrac12$, so lowering the base point from $\alpha'$ to $\alpha$ deconvolves by a Cauchy kernel of width $\alpha' - \alpha$; the Cauchy kernels form a semigroup under convolution, and the descent $\alpha\downarrow\tfrac12$ is a backward Cauchy flow. Under RH it ends at the atoms $\delta_{\gamma^2}$, and by Theorem~\ref{thm:continuity} the laws $\mu^{(\alpha)}$ converge weakly to $\mu$ as their Thorin measures collapse. Table~\ref{tab:descent} and Figure~\ref{fig:descent} show this on the first zeros.

\begin{table}[htbp]
\centering
\small
\begin{tabular}{rrrrr}
\toprule
$\alpha$ & $T = 20$ & $T = 30$ & $T = 40$ & $T = 50$\\
\midrule
$1$ & $1.2377$ & $3.4003$ & $6.2304$ & $9.7066$\\
$0.75$ & $1.1228$ & $3.2317$ & $6.1202$ & $9.7694$\\
$0.6$ & $1.0496$ & $3.0986$ & $6.0488$ & $9.8816$\\
$0.52$ & $1.0099$ & $3.0200$ & $6.0098$ & $9.9747$\\
$N(T)$ & $1$ & $3$ & $6$ & $10$\\
\bottomrule
\end{tabular}
\caption{The Thorin distribution $U_\alpha((0,T^2])$ of \eqref{eq:arctan} converging to the zero count as the base point descends. At $\alpha = 0.75$, $T = 50$ the arctangent sum gives $9.7694$ and the argument of $\xi$ along the line gives $9.7775$, the difference being the zeros beyond the $600$th.}
\label{tab:descent}
\end{table}

An off line zero breaks the flow. If $\beta > \alpha$ then $c_\rho < 0$, its arctangent pair decreases by $\pi$ instead of increasing, and the zero is counted with the wrong sign; the measure stops being positive at exactly $\alpha = \beta$. So the descent is possible down to $\sigma_{\max}$ and no further, and
\[
\sigma_{\max} - \tfrac12 = \inf\{\alpha - \tfrac12 : U_\alpha \ge 0\}
\]
is a Cauchy analogue of the de Bruijn Newman constant. The two are related: if all zeros have displacement at most $\delta_{\max} = \sigma_{\max} - \tfrac12$ then $\Lambda \le 2\delta_{\max}^2$ in the normalisation of \cite{Ladder}, by de Bruijn's theorem~\cite{deBruijn}. The Gaussian flow sees the displacement quadratically, the Cauchy flow linearly, and both vanish exactly under RH.

\subsection*{What survives below \texorpdfstring{$\alpha = 1$}{alpha = 1}}

Extending the laws $\mu^{(\alpha)}$ as GGCs to some fixed $\alpha < 1$ is equivalent to the quasi Riemann hypothesis $\sigma_{\max} \le \alpha$, which is open for every $\alpha < 1$. What can be said unconditionally is a statement up to height $H$. For $\alpha \in (\tfrac12, 1)$ split the product of Theorem~\ref{thm:product} as $\varphi^{(\alpha)} = \varphi^{(\alpha)}_{\le H}\,R^{(\alpha)}$, where the first factor runs over $\gamma \le H$. Every zero below $H$ lies on the line, so $c_\rho = \alpha - \tfrac12 > 0$ and $\varphi^{(\alpha)}_{\le H}$ is a GGC Laplace transform with Thorin distribution given by \eqref{eq:arctan} restricted to $\gamma \le H$. The remainder satisfies, for real $s \ge 0$,
\[
\big|\log R^{(\alpha)}(s)\big| \le \sum_{\gamma>H}\frac{2|c_\rho|\sqrt s + s}{\gamma^2}\,(1+o(1)) \le \big(\sqrt s + s\big)\frac{\log H}{2\pi H}(1+o(1)),
\]
using $|c_\rho| \le \tfrac12$ and the zero counting function. So on $[0, s_0]$ the reciprocal at any base point in the strip is within $(\sqrt{s_0}+s_0)\log H/(2\pi H) \approx 1.5\times10^{-12}(\sqrt{s_0}+s_0)$ in log Laplace transform of an explicit GGC. By Theorem~\ref{thm:continuity} this is not enough: closure needs exact limits, and $H \to \infty$ is not available.

\subsection*{A third description: the argument of \texorpdfstring{$\xi$}{xi} as a clock}

Formula \eqref{eq:arctan} has a consequence that turns the descent into a statement about a single analytic function on vertical lines. Differentiating the argument of $\xi$ along $\Real z = \alpha$ gives $\frac{d}{d\tau}\arg\xi(\alpha+i\tau) = \Imag\big(i\,\xi'/\xi\big) = \Real\,\xi'/\xi(\alpha+i\tau)$, so in the ordinate variable $\tau = \sqrt z$ the Thorin density of $\mu^{(\alpha)}$ is
\begin{equation}\label{eq:lagdensity}
\frac{1}{\pi}\Real\frac{\xi'}{\xi}(\alpha+i\tau) = \frac1\pi\Big[\Real\Big(\frac1s+\frac{1}{s-1}\Big) + \frac12\Real\frac{\Gamma'}{\Gamma}\Big(\frac s2\Big) - \frac12\log\pi + \Real\frac{\zeta'}{\zeta}(s)\Big], \qquad s = \alpha+i\tau .
\end{equation}
This identity is unconditional: the right side is the signed density defined by the argument, whatever the zeros do. It gives a representation of the base point law that is positive by construction whenever the argument is monotone. Every GGC with Thorin distribution function $F_U(t) = U((0,t])$ is the law of $\int_0^\infty dg_v/F_U^{-1}(v)$, where $g$ is the standard gamma subordinator, with Lévy density $e^{-x}/x$, because the Laplace exponent of $\int h\,dg$ is $\int\log(1+sh(v))\,dv$. With $F_U$ given by \eqref{eq:arctan},
\begin{equation}\label{eq:clock}
\mu^{(\alpha)} = \text{law of }\int_0^\infty\frac{dg_v}{\tau_\alpha(v)^2}, \qquad \tau_\alpha = \text{the inverse of } \tau \mapsto \frac1\pi\Delta_0^\tau\arg\xi(\alpha+i\tau).
\end{equation}
The argument of $\xi$ along the line is a clock, and $\mu^{(\alpha)}$ is a gamma process run against that clock. At $\alpha = \tfrac12$ under RH the clock is the step function $N(\tau)$, its inverse takes the value $\gamma_n$ on $(n-1,n]$, the gamma increments over unit intervals are independent standard exponentials, and \eqref{eq:clock} becomes $\sum E_n/\gamma_n^2$. The representation exists exactly when the clock runs forward, that is when \eqref{eq:lagdensity} is nonnegative on the whole line.

\begin{theorem}\label{thm:lagarias}
For $\alpha_0 \ge \tfrac12$ the following are equivalent.
\begin{enumerate}
\item[(i)] No zero of $\zeta$ has real part greater than $\alpha_0$.
\item[(ii)] $\mu^{(\alpha)}$ is a GGC for every $\alpha > \alpha_0$.
\item[(iii)] The clock $\tau\mapsto\arg\xi(\alpha+i\tau)$ is nondecreasing on $[0,\infty)$ for every $\alpha > \alpha_0$, so that \eqref{eq:clock} holds.
\item[(iv)] $\Real\,\xi'/\xi(s) > 0$ for every $s$ with $\Real s > \alpha_0$.
\end{enumerate}
For $\alpha_0 = \tfrac12$ condition (iv) is Lagarias' criterion for RH \cite{Lagarias}, and all four are equivalent to RH.
\end{theorem}

\begin{proof}
(i) implies (ii) is Theorem~\ref{thm:base}, and (ii) implies (iii) because a Thorin measure is nonnegative and its distribution function is the clock, by \eqref{eq:arctan}. (iii) implies $\Real\,\xi'/\xi \ge 0$ on $\Real s > \alpha_0$, by \eqref{eq:lagdensity} and the symmetry $\xi(\overline s) = \overline{\xi(s)}$. This weak inequality implies (i): if a zero $\rho = \beta + i\gamma$ of multiplicity $m$ had $\beta > \alpha_0$, then at $s = \alpha + i\gamma$ with $\alpha_0 < \alpha < \beta$ the Hadamard sum
\[
\Real\frac{\xi'}{\xi}(s) = \sum_\rho\frac{\sigma - \beta_\rho}{|s-\rho|^2}
\]
contains the term $-m/(\beta-\alpha)$, while the remaining terms stay bounded as $\alpha \uparrow \beta$, so the sum is negative for $\alpha$ close to $\beta$. Finally (i) implies (iv), since under (i) every term of the Hadamard sum is positive when $\sigma > \alpha_0$. So (i) to (iii) are equivalent to each other, (i) implies (iv), and (iv) implies the weak inequality, hence (i).
\end{proof}

So a third description of $U_\alpha$, beyond the prime side and the zero side, exists, and it is the most classical one: the argument of $\xi$ along vertical lines. What it shows is that extending the laws into the strip is the same problem as Lagarias' positivity of $\Real\,\xi'/\xi$ on half planes, and the decomposition \eqref{eq:lagdensity} says where the difficulty lies. The archimedean part is $\tfrac12\log(\tau/2\pi) + O(1/\tau)$ for large $\tau$; the clock runs forward at height $\tau$ exactly when
\begin{equation}\label{eq:onesided}
-\Real\frac{\zeta'}{\zeta}(\alpha+i\tau) < \frac12\log\frac{\tau}{2\pi} + O\Big(\frac1\tau\Big),
\end{equation}
a one sided bound on the logarithmic derivative of $\zeta$. The known regions where it holds unconditionally are these.
\begin{enumerate}
\item[(a)] For $\alpha > 1$, $-\Real\,\zeta'/\zeta(\alpha+i\tau) \le -\zeta'/\zeta(\alpha) = \sum\Lambda(n)n^{-\alpha}$, which is about $1/(\alpha-1)$, so the Euler product alone gives \eqref{eq:onesided} for $\alpha \ge 1 + 2/\log(\tau/2\pi) + O(\cdot)$, a region that approaches the line $\Real s = 1$ from the right as $\tau \to \infty$.
\item[(b)] For $\alpha \ge 1$ at every height, by Theorem~\ref{thm:base}, since no zero has $\beta \ge 1$; this closes the gap left by (a).
\item[(c)] Below $1$, in the Vinogradov Korobov~\cite{Ford} zero free region $\alpha \ge 1 - c\,(\log\tau)^{-2/3}(\log\log\tau)^{-1/3}$, where $\zeta'/\zeta(\alpha+i\tau) \ll (\log\tau)^{2/3}(\log\log\tau)^{1/3} = o(\log\tau)$, so \eqref{eq:onesided} holds at every sufficiently large height. This region approaches the line $\Real s = 1$ from the left.
\end{enumerate}
The picture is symmetric about $\Real s = 1$: the Euler product controls the clock to the right of a curve converging to the line from the right, and the zero free region controls it to the left of a curve converging from the left. At no fixed $\alpha < 1$ is the clock known to run forward at all heights, and making it do so is the quasi Riemann hypothesis.

Numerically, on the lines $\alpha = 1, 0.75, 0.6, 0.55, 0.51$ and heights up to $200$ to $400$, $\Real\,\xi'/\xi$, which is $\pi$ times the density \eqref{eq:lagdensity}, is positive, with minima $0.0232, 0.0116, 0.0053, 0.0027, 0.0005$ attained at the lowest heights, below the first zero. There it is close to its value $\mathsf B_\alpha = \xi'(\alpha)/\xi(\alpha)$ at $\tau = 0$, and it is a cancellation: at $\alpha = 0.55$, $\tau = 5.1$ the archimedean part is $-0.102$ and the zeta part $+0.105$. As a check of \eqref{eq:lagdensity} and the Thorin representation, integrating $\log(1+s/t)$ against the density in $t = \tau^2$ at $\alpha = 1$, $s = 1$ gives $0.046125$ against $\log\xi(2)/\xi(1) = 0.046118$, the difference being the asymptotic tail used beyond $t = 10^4$.

\subsection*{The arithmetic content, and the sanity tests}

The statement that $\mu^{(\alpha)}$ is a GGC for $\alpha \ge 1$ uses the absence of zeros with $\Real\rho > 1$, which comes from the Euler product, and the absolute continuity of its Thorin measure at $\alpha = 1$ uses the absence of zeros on $\Real\rho = 1$, which is the prime number theorem (Section~\ref{sec:alpha1}). Both fail for the Davenport Heilbronn function~\cite{DavenportHeilbronn}, which has zeros with real part above $1$, and both hold, with $1$ replaced by $\tfrac12 + \delta_{\max}(t)$, for the heat flowed $H_t$. So Theorem~\ref{thm:base} at $\alpha \ge 1$ is the one positivity statement in these notes that genuinely uses arithmetic. RH asks for the same statement at every $\alpha > \tfrac12$, and by the closure of Theorem~\ref{thm:continuity} it is enough to have it on a sequence $\alpha_n \downarrow \tfrac12$.

What would extend the laws into the strip is therefore an arithmetic proof that $U_\alpha \ge 0$ for some $\alpha < 1$, a quasi RH, by a route other than locating zeros. From the prime side, $U_\alpha$ is the sine squared transform of the signed Lévy measure of Section~\ref{sec:sin2}, which converges only for $\alpha > 1$; below $1$ the prime series must be continued analytically, and analytic continuation preserves identities but not positivity. From the zero side, $U_\alpha$ is a backward Cauchy flow of the counting measure, which is positive only if the zeros are where RH says they are. A third description exists, the gamma process clock of \eqref{eq:clock}, but it is positive only where the argument of $\xi$ is monotone on vertical lines, which Theorem~\ref{thm:lagarias} shows is the same condition. A proof would need a reason for the one sided bound \eqref{eq:onesided} at some fixed $\alpha < 1$ and all heights, coming from neither the convergent Euler product nor the location of zeros.

\section{The base point \texorpdfstring{$\alpha = 1$}{alpha = 1}}\label{sec:alpha1}

The line $\Real s = 1$ is where the arithmetic of $\zeta$ meets the base point picture, and the law at $\alpha = 1$ deserves to be written out on its own. Since $\xi(1) = \tfrac12$, its Laplace transform is simply
\[
\varphi^{(1)}(s) = \frac{\xi(1)}{\xi(1+\sqrt s)} = \frac{1}{2\,\xi(1+\sqrt s)} .
\]

\begin{theorem}[The law at $\alpha = 1$]\label{thm:alpha1}
Unconditionally, $\varphi^{(1)}$ is the Laplace transform of a GGC law $\mu^{(1)}$ with zero drift, and the following hold.
\begin{enumerate}
\item[(i)] Its Thorin measure is $U_1 = \sum_{\gamma>0}u_{\rho,1}(z)\,dz$ with $c_\rho = 1 - \beta_\rho \in (0,1)$ in \eqref{eq:lorentz}; in the ordinate variable $\tau = \sqrt z$ its density is $\pi^{-1}\Real\,\xi'/\xi(1+i\tau)$, and its distribution function is $\pi^{-1}\Delta_0^\tau\arg\xi(1+i\tau)$.
\item[(ii)] Its Lévy measure is $k_1(x)x^{-1}dx$ with $k_1(x) = \sum_{\gamma>0}\Real\big[e^{(1-\beta+i\gamma)^2x}\operatorname{erfc}((1-\beta+i\gamma)\sqrt x)\big]$.
\item[(iii)] It has a stable tail of index one half, $\mu^{(1)}((y,\infty)) \sim \mathsf B_1/\sqrt{\pi y}$, with the explicit constant
\[
\mathsf B_1 = \frac{\xi'(1)}{\xi(1)} = \sum_\rho\frac1\rho = 1 + \frac{\gamma_{\rm E}}{2} - \frac12\log4\pi = 0.0230957\ldots,
\]
so $\mu^{(1)}$ has infinite mean.
\item[(iv)] The Thorin measure $U_1$ has no atoms, that is, it is absolutely continuous, if and only if $\zeta(1+i\tau) \ne 0$ for all real $\tau$.
\end{enumerate}
\end{theorem}

\begin{proof}
The GGC property and (i) to (iii) are the case $\alpha = 1$ of Theorem~\ref{thm:base}, \eqref{eq:erfc} and the tail computation that follows it. Only the absence of zeros with $\Real\rho > 1$ is used for the GGC property, and that comes from the absolute convergence of the Euler product. The value of $\mathsf B_1$ is the classical evaluation of $\sum_\rho 1/\rho$, summed symmetrically, which equals $\xi'(1)/\xi(1) = -\xi'(0)/\xi(0)$ by the functional equation. For (iv), the atoms of $U_\alpha$ come exactly from zeros with $c_\rho = 0$, that is with $\beta_\rho = \alpha$: such a zero contributes the factor $\gamma^2/(s+\gamma^2)$, whose Thorin measure is $\delta_{\gamma^2}$, while every zero with $c_\rho > 0$ contributes the absolutely continuous density \eqref{eq:lorentz}.
\end{proof}

Part (iv) is the prime number theorem in Thorin form. The GGC property at $\alpha = 1$ is cheap: it needs only the Euler product in the half plane where it converges absolutely, and it would hold equally if $\zeta$ had zeros on the line $\Real s = 1$. What the theorem of Hadamard and de la Vallée Poussin adds is that those zeros do not exist, and in Thorin terms this is the statement that the Thorin measure at $\alpha = 1$ is spread out by Cauchy kernels of strictly positive width rather than concentrated in atoms. Equivalently, the clock $\tau\mapsto\arg\xi(1+i\tau)$ is continuous: it never jumps. The prime number theorem is thus the first step of the descent of Section~\ref{sec:base}, the step that moves every Cauchy width off zero at $\alpha = 1$; RH is the statement that the widths $\alpha - \beta_\rho$ stay nonnegative all the way down to $\alpha = \tfrac12$, where they all vanish together and the measure becomes the atoms $\delta_{\gamma^2}$.

Under RH every width at $\alpha = 1$ equals $\tfrac12$, and $U_1$ is the zero counting measure seen through a Cauchy window of width $\tfrac12$ in the ordinate:
\[
\frac{1}{\pi}\Real\frac{\xi'}{\xi}(1+i\tau) = \frac1\pi\sum_{\gamma>0}\Big[\frac{1/2}{\tfrac14+(\tau-\gamma)^2} + \frac{1/2}{\tfrac14+(\tau+\gamma)^2}\Big] \quad\text{under RH}.
\]
Table~\ref{tab:alpha1} lists the density. It starts at $\mathsf B_1/\pi = 0.00735$ at $\tau = 0$, stays near that value below the first zero, rises to $0.647$ at $\tau = \gamma_1 = 14.13$, where the window sits on a zero, and thereafter fluctuates around the smooth value $(2\pi)^{-1}\log(\tau/2\pi)$, since the zero spacing at these heights, about $2.5$, is several window widths. As a check of the whole representation, integrating $\log(1+s/\tau^2)$ against this density at $s = 1$ gives $0.046125$, against the exact value $\log(\xi(2)/\xi(1)) = 0.046118$; the difference is the asymptotic tail used beyond $\tau = 100$. The product formula of Theorem~\ref{thm:product} at $\alpha = 1$ reproduces $\xi(1)/\xi(1+w)$ to six digits at $w = 0.5, 1, 3$, for instance $0.954929659$ against $0.954929048$ at $w = 1$.

\begin{table}[htbp]
\centering
\small
\setlength{\tabcolsep}{4pt}\begin{tabular}{rrrrrrrrrr}
\toprule
$\tau$ & $0$ & $5$ & $10$ & $14.13$ & $20$ & $30$ & $50$ & $100$ & $200$\\
\midrule
$\pi^{-1}\Real\,\xi'/\xi(1+i\tau)$ & $0.0074$ & $0.0084$ & $0.0168$ & $0.6474$ & $0.1419$ & $0.4071$ & $0.6050$ & $0.2308$ & $0.2134$\\
$(2\pi)^{-1}\log(\tau/2\pi)$ & & & $0.074$ & $0.129$ & $0.184$ & $0.249$ & $0.330$ & $0.440$ & $0.551$\\
\bottomrule
\end{tabular}
\caption{The Thorin density of $\mu^{(1)}$ in the ordinate variable, against the smooth zero density.}
\label{tab:alpha1}
\end{table}

\begin{figure}[htbp]
\centering
\includegraphics[width=0.88\textwidth]{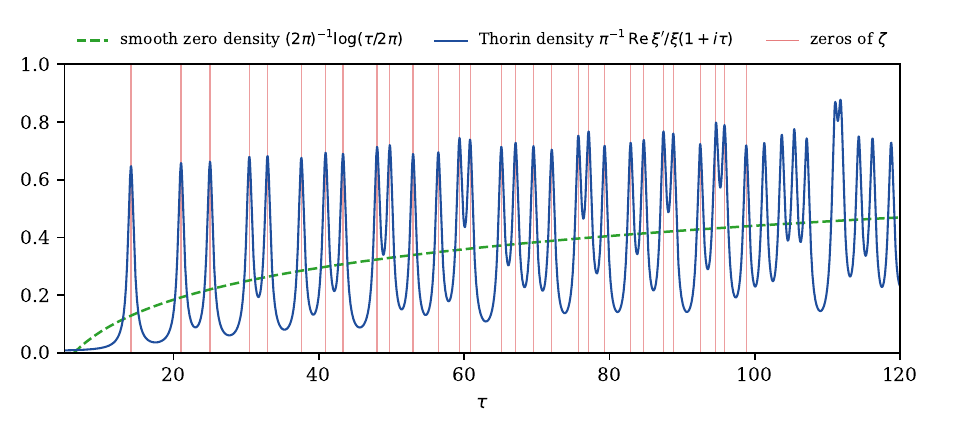}
\caption{The Thorin density of the law at $\alpha = 1$ in the ordinate variable (blue), with the smooth zero density $(2\pi)^{-1}\log(\tau/2\pi)$ (green, dashed) and the zeros of $\zeta$ (red vertical lines). Under RH the density is the zero count seen through a Cauchy window of width one half.}
\label{fig:alpha1}
\end{figure}

Finally, the one sided bound \eqref{eq:onesided} at $\alpha = 1$ reads $-\Real\,\zeta'/\zeta(1+i\tau) < \tfrac12\log(\tau/2\pi) + O(1/\tau)$. It holds at every height, because $U_1 \ge 0$; at large heights it also follows directly from the Vinogradov Korobov bound $\zeta'/\zeta(1+i\tau) \ll (\log\tau)^{2/3}(\log\log\tau)^{1/3}$, which leaves a margin of order $\log\tau$. That margin is what makes the strip $1 - c(\log\tau)^{-2/3}(\log\log\tau)^{-1/3} \le \alpha \le 1$ accessible at large heights, and its disappearance as $\tau\to\infty$ at any fixed $\alpha < 1$ is the quasi Riemann hypothesis.

\subsection*{The prime side at \texorpdfstring{$\alpha = 1$}{alpha = 1}}

Problem 2 of Section~\ref{sec:open} asks for a proof of the positivity of $U_1$ from the primes rather than from the zeros. This subsection records what can be said about it: an unconditional Poisson structure that ties $\alpha = 1$ to all $\alpha > 1$, an exact prime side form of the inequality, the region where a one line prime side proof exists, and numerical evidence for the margin.

\begin{theorem}[Poisson structure above the line]\label{thm:poisson}
Extend $U_\alpha$ to an even measure on $\R$ in the ordinate variable, $U_\alpha(d\tau) = \pi^{-1}\Real\,\xi'/\xi(\alpha+i|\tau|)\,d\tau$. Then, unconditionally, for every $\alpha > 1$,
\[
U_\alpha = P_{\alpha-1} * U_1, \qquad P_c(x) = \frac{c}{\pi(c^2+x^2)},
\]
so every base point law to the right of $1$ is the Cauchy smoothing of the law at $1$. Equivalently, $s \mapsto \xi'/\xi(1+s)$ has positive real part on $\Real s > 0$ and is the Poisson integral of its boundary values, with $\pi U_1$ as boundary measure and no singular part.
\end{theorem}

\begin{proof}
By \eqref{eq:arctan} the density of $U_\alpha$ is $\sum_\rho[P_{c_\rho(\alpha)}(\tau-\gamma) + P_{c_\rho(\alpha)}(\tau+\gamma)]$ with $c_\rho(\alpha) = \alpha - \beta_\rho = c_\rho(1) + (\alpha-1)$, and $c_\rho(1) = 1-\beta_\rho > 0$ for every zero. The Cauchy kernels form a convolution semigroup, $P_a * P_b = P_{a+b}$, so each term at $\alpha$ is the convolution of the corresponding term at $1$ with $P_{\alpha-1}$, and the sums converge by Theorem~\ref{thm:product}. The second statement is the same fact read as the Herglotz representation of a function with positive real part on a half plane; the absence of a singular boundary part is the absence of zeros on $\Real s = 1$.
\end{proof}

So the descent from any $\alpha > 1$ to $\alpha = 1$ is harmless: it is an inverse Poisson problem whose solution is known to exist and to be positive. The whole difficulty of the descent is concentrated at and below $1$.

The inequality itself has an exact prime side form. On the line $\Real s = 1$ the archimedean part of \eqref{eq:lagdensity} is $\tfrac12\log(\tau/2\pi) + O(1/\tau)$, and by Mertens' theorem in the form that follows from the prime number theorem, $\sum_{n\le X}\Lambda(n)/n = \log X - \gamma_{\rm E} + o(1)$. Taking $X = \tau/2\pi$,
\[
\text{archimedean part at } 1+i\tau = \frac12\sum_{n\le\tau/2\pi}\frac{\Lambda(n)}{n} + \frac{\gamma_{\rm E}}{2} + o(1),
\]
which the table below confirms to two decimals. Positivity of $U_1$ is therefore the statement
\begin{equation}\label{eq:mertensform}
\sum_{n\ge2}\frac{\Lambda(n)}{n}\cos(\tau\log n) \;\le\; \frac12\sum_{n\le\tau/2\pi}\frac{\Lambda(n)}{n} + \frac{\gamma_{\rm E}}{2} + o(1),
\end{equation}
the left side being the conditionally convergent series for $-\Real\,\zeta'/\zeta(1+i\tau)$. In words: the oscillating prime sum at height $\tau$ is at most half the Mertens sum up to the cutoff $\tau/2\pi$. The cutoff is the Riemann Siegel scale: for $n \le \tau/2\pi$ consecutive phases $\tau\log n$ differ by $\tau/n \ge 2\pi$, so the terms are arithmetically scrambled, while for $n > \tau/2\pi$ the phase moves slowly and the tail cancels by smooth summation. The inequality has a great deal of room. If the phases on $n\le\tau/2\pi$ were independent, the left side would have mean zero and variance $\tfrac12\sum\Lambda(n)^2/n^2 < \infty$, so it would be $O(1)$ against a right side of size $\tfrac12\log\tau$; under RH the left side is in fact $O(\log\log\tau)$. The right side is attained only if the phases conspire over a long range of $n$, and such a conspiracy is exactly a zero of $\zeta$ near $1+i\tau$. That is why the zero side proof is one line and the prime side proof is not.

Where does a prime side argument reach today? For $\alpha > 1$ the trivial bound $\sum\Lambda(n)n^{-\alpha}\cos(\tau\log n) \le -\zeta'/\zeta(\alpha) \le 1/(\alpha-1)$ proves positivity of $U_\alpha$ at height $\tau$ whenever $\alpha - 1 \ge 2/\log(\tau/2\pi)$, with no information about zeros at all; by Theorem~\ref{thm:poisson} this is positivity of the Cauchy smoothings of $U_1$ of width at least $2/\log(\tau/2\pi)$. Smoothing in the ordinate rather than by moving the base point does better, and the reason is a small miracle of bookkeeping: a smoothing kernel of bandwidth $L$ truncates the prime sum at $n \le e^L$, because the $n$th term oscillates in $\tau$ at frequency $\log n$.

\begin{theorem}[Prime side positivity at the zero spacing scale]\label{thm:fejer}
Let $K \ge 0$ be an even probability density on $\R$ whose Fourier transform $\widehat K(\lambda) = \int K(x)e^{i\lambda x}dx$ vanishes for $|\lambda| > L$. Then, in the ordinate variable,
\[
\pi\,(K * U_1)(\tau) = (K*\mathrm{arch})(\tau) - \sum_{n\le e^L}\frac{\Lambda(n)}{n}\,\widehat K(\log n)\cos(\tau\log n), \qquad \mathrm{arch}(\tau) = \tfrac12\log\frac{\tau}{2\pi} + O\Big(\frac1{\tau}\Big),
\]
and $(K*\mathrm{arch})(\tau) = \mathrm{arch}(\tau) + O(\log\tau/(L\tau))$. Consequently $(K*U_1)(\tau) \ge 0$ whenever
\[
\sum_{n\le e^L}\frac{\Lambda(n)}{n}\,|\widehat K(\log n)| \le \frac12\log\frac{\tau}{2\pi} - O\Big(\frac{\log\tau}{L\tau}\Big),
\]
an inequality that involves only the primes. For the Fejér kernel, $\widehat K(\lambda) = (1-|\lambda|/L)_+$, Mertens' theorem gives $\sum_{n\le e^L}\Lambda(n)n^{-1}(1-\log n/L) = \tfrac12L - \gamma_{\rm E} + o(1)$, so
\[
(F_L * U_1)(\tau) \ge 0 \qquad\text{for all } L \le \log\frac{\tau}{2\pi} + 2\gamma_{\rm E} - o(1),
\]
and at the threshold $L = \log(\tau/2\pi)$ the main lobe of $F_L$ extends to $\pm2\pi/L$, which is the mean spacing $2\pi/\log(\tau/2\pi)$ of the zeros at height $\tau$. For the modulated kernel $\tfrac12[F_L(x-2\pi/L) + F_L(x+2\pi/L)]$, with $\widehat K(\lambda) = (1-|\lambda|/L)_+\cos(2\pi\lambda/L)$, the same argument gives $L \le \tfrac\pi2\log(\tau/2\pi) + O(1)$.
\end{theorem}

\begin{proof}
For $\sigma > 1$ the series $-\zeta'/\zeta(\sigma+i\tau) = \sum\Lambda(n)n^{-\sigma-i\tau}$ converges absolutely and may be convolved term by term; since $K$ is even, convolving $\cos(\tau\log n)$ with $K$ gives $\widehat K(\log n)\cos(\tau\log n)$, which vanishes for $n > e^L$, and $n^{-\sigma} \le n^{-1}$. The passage to $\sigma = 1$ uses the existence and continuity of $\zeta'/\zeta$ on $\Real s = 1$, $\tau \ne 0$, which is the prime number theorem, together with the classical bound $\zeta'/\zeta(\sigma+it) \ll \log|t|$ for $\sigma \ge 1$ to dominate the tails of $K$; the inequality itself uses nothing about zeros. The archimedean part is smooth on the scale of $K$, which has width $1/L$ and tails $O(1/(Lu^2))$, and is bounded by $O(\log\tau)$ on the whole line, so $K*\mathrm{arch} = \mathrm{arch} + O(\log\tau/(L\tau))$. Bounding each cosine by $1$ gives the criterion. For the Fejér weights, partial summation against $\sum_{n\le X}\Lambda(n)/n = \log X - \gamma_{\rm E} + o(1)$ gives $\int_0^L(1-\lambda/L)\,d(\lambda - \gamma_{\rm E}) = \tfrac12L - \gamma_{\rm E}$ up to $o(1)$; for the modulated weights $\int_0^L(1-\lambda/L)|\cos(2\pi\lambda/L)|\,d\lambda = L/\pi$. With only the elementary Mertens bound $\sum_{n\le X}\Lambda(n)/n = \log X + O(1)$ the thresholds hold with $2\gamma_{\rm E}$ replaced by an absolute constant.
\end{proof}

So the prime side, using Mertens' theorem and the prime number theorem but no further information about the zeros, proves that the Thorin measure at $\alpha = 1$ is positive on every window of about one mean zero spacing, and the constants can be pushed by a bounded factor but not beyond order $\log\tau$. This is the natural limit of a trivial bound: a window narrower than the mean spacing resolves individual zeros, and whether the measure is positive there is exactly whether the zero near the window sits on the line or not. Beyond the threshold the inequality must gain over the trivial bound by the factor $L/\log\tau$, and all known ways to gain it go through the zeros. Problem 2 is therefore the problem of making the prime sum $\sum_{n\le e^L}\Lambda(n)n^{-1}\widehat K(\log n)\cos(\tau\log n)$ cancel, on average over $n$, once $L$ exceeds $\log\tau$ by more than a constant factor; whatever argument does this at $\alpha = 1$ is the candidate for problem 1 one step to the left, where the same sum appears with $n^{-\alpha}$ in place of $n^{-1}$ and the archimedean term unchanged.

Numerically the margin is comfortable, and the head of the series does the work. On $2 \le \tau \le 2000$ the left side of \eqref{eq:mertensform}, computed as $-\Real\,\zeta'/\zeta(1+i\tau)$, has maximum $1.588$ at $\tau = 1967$, where the right side is $2.873$; its largest ratio to the right side above $\tau = 20$ is $0.744$, at $\tau = 27.8$; and its mean over $[100,2000]$ is $0.000$, against a mean right side of $2.46$. The sharp partial sum $\sum_{n\le\tau/2\pi}\Lambda(n)n^{-1}\cos(\tau\log n)$ tracks the full value closely, with correlation $0.974$ on $[20,3000]$ and a difference of standard deviation $0.14$ against $0.63$ for the sum itself, so the fluctuation lives below the Riemann Siegel cutoff. This is what the truncated explicit formula predicts: $\sum_{n\le x}\Lambda(n)n^{-s} = -\zeta'/\zeta(s) + x^{1-s}/(1-s) - \sum_\rho x^{\rho-s}/(\rho-s) + \ldots$, and at $x = \tau/2\pi$, $s = 1+i\tau$ the zero sum is $O(x^{-1/2}\log^2\tau)$ under RH. The Fejér weighted head at $L = \log(\tau/2\pi)$ has maximum $1.33$ on the same range, against a trivial bound $\tfrac12L - \gamma_{\rm E}$ between $2.3$ and $2.5$ and an archimedean term between $2.9$ and $3.1$ at the top of the range, confirming the theorem with room to spare.

\begin{table}[htbp]
\centering
\small
\begin{tabular}{rrrr}
\toprule
$\tau$ & archimedean part & $\tfrac12\sum_{n\le\tau/2\pi}\Lambda(n)/n$ & $\tfrac12\sum + \gamma_{\rm E}/2$\\
\midrule
$50$ & $1.037$ & $0.743$ & $1.032$\\
$200$ & $1.730$ & $1.472$ & $1.761$\\
$1000$ & $2.535$ & $2.244$ & $2.533$\\
$2000$ & $2.882$ & $2.602$ & $2.891$\\
\bottomrule
\end{tabular}
\caption{The archimedean part of the Thorin density at $\alpha = 1$ against half the Mertens sum to the Riemann Siegel cutoff.}
\end{table}

\begin{figure}[htbp]
\centering
\includegraphics[width=0.85\textwidth]{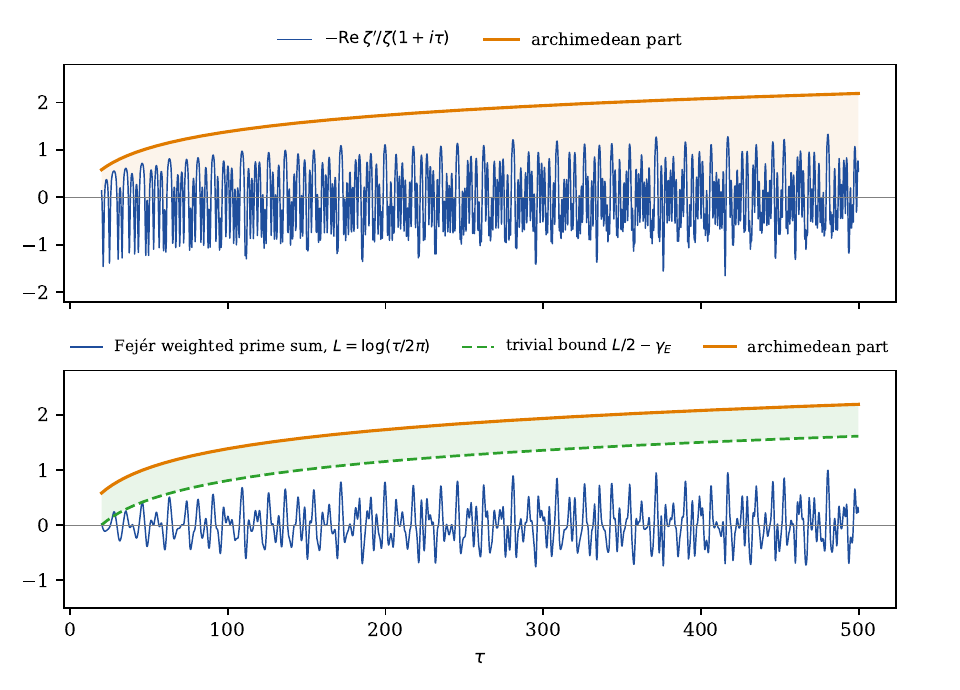}
\caption{Problem 2 in pictures. Above: the two sides of the Mertens form \eqref{eq:mertensform}, the oscillating prime sum $-\Real\,\zeta'/\zeta(1+i\tau)$ (blue) and the archimedean part (orange); positivity of $U_1$ is the blue curve staying below the orange one, and the shaded region is the margin. Below: the Fejér weighted prime sum at bandwidth $L = \log(\tau/2\pi)$ (blue), the trivial bound $L/2 - \gamma_{\rm E}$ of Theorem~\ref{thm:fejer} (green, dashed) and the archimedean part (orange); the green band is the margin by which the theorem holds, and the gap between the blue curve and the green line is what the trivial bound gives away.}
\label{fig:mertens}
\end{figure}

\section{Is there a path to RH?}\label{sec:path}

The honest answer is: not yet, but the question has acquired a precise shape, and the shape tells us what any path must look like.

The situation is this. The reciprocal Xi law $\mu$ is unconditionally self decomposable, and its canonical function $\Theta$ is unconditionally monotone to order $H^2/80 \approx 10^{23}$ at every resolution and to every order at resolution $x \le 1/20$. RH is the statement that $\Theta$ is monotone to every order at every resolution, that is, completely monotone, that is, the Laplace transform of a positive Thorin measure. Between what is known and what is wanted lies a single limit, $k \to \infty$ at fixed $x > 1/20$. Everything proved so far stops short of that limit for one reason: the proofs use the zeros (verified up to $H$) and the counting function, together with the fact that any zero has $|\delta| < 1/2$. They do not use the arithmetic beyond bounding the size of the prime sum near $x = 0$. Any argument of that kind is bound to stop at a finite order, because the same arguments apply to functions that have zeros off the line: to the heat flowed $H_t$ for $t < 0$, whose zeros are all real up to a great height and then collide, and to the Davenport Heilbronn function, which has the functional equation and the gamma factor of $\xi$ but no Euler product. For both of these $\Theta$ is positive, self decomposable, monotone to a high order, and not completely monotone. A proof of RH must therefore use something that these two functions lack: the Euler product at every order, not merely in the small $x$ regime where the primes are exponentially suppressed.

Where does the Euler product enter the present framework? In exactly two places. First, through the positivity $\Lambda(n) \ge 0$ of the prime measure, which makes the prime part of the Lévy density a negative superposition of first passage densities, equation \eqref{eq:levydensity}, and which, after an Esscher tilt, makes the zeta distributed level $\log N$ a genuine random variable, so that the balance
\[
Y_\lambda + T_{\log N} + T_\kappa \stackrel{d}{=} E_{\lambda-1/4} + C_\lambda
\]
of \cite{Ladder} holds, with every term except $Y_\lambda$ an explicit Brownian functional. Second, through the sine squared transform at base points $\alpha > 1$, where the Euler product converges and the prime part of the Thorin density is manifestly positive, so that the reciprocal $\xi(\alpha)/\xi(\alpha+\sqrt s)$ is a GGC after a drift correction and RH is the survival of that positive Thorin measure down to $\alpha = 1/2$. These are the only two structures in the whole picture that fail for the heat flow and for Davenport Heilbronn, and so they are the only candidates for the arithmetic input a proof needs.

Three routes present themselves, and each has a clear obstacle.

The first route is the balance. RH is the statement that the deconvolution $Y_\lambda = (E + C_\lambda) \ominus (T_{\log N} + T_\kappa)$ is a GGC, equivalently that the canonical function $e^{-\lambda x}(A - P)$ is completely monotone, where $A$ is archimedean and $P$ is the prime sum. The archimedean side has an explicit signed Thorin measure, with a unit atom at $\lambda - 1/4$ and negative mass $-1.124$ just above $\lambda$; the prime side is a positive combination of first passage densities whose Thorin side is not a measure. So the statement to prove is that a positive Lévy measure built from first passage times to the levels $\log n$, with weights $\Lambda(n)/(\sqrt n\log n)$, exactly removes the negative Thorin mass of the gamma factor and introduces no negative mass of its own. This is a statement about primes, with the gamma factor as the known side, and it is equivalent to RH. The obstacle is that complete monotonicity is a condition on all derivatives, and the prime side at order $k$ and resolution $x$ is a Gaussian smoothed Dirichlet polynomial at height $\sqrt{k/x}$, which is to say the fluctuation $S(T)$ of the zero count; controlling it at every order is controlling the zeros.

The second route is the base point descent of Section~\ref{sec:base}. The plain reciprocal $\xi(\alpha)/\xi(\alpha+\sqrt s)$ is a GGC for every $\alpha \ge 1$, unconditionally, with Thorin measure the zero counting measure smoothed by Cauchy kernels of width $\alpha - \beta_\rho$, and RH is the statement that this backward Cauchy flow stays positive down to $\alpha = \tfrac12$; by Bondesson's continuity theorem it is enough to reach a sequence $\alpha_n \downarrow \tfrac12$. Section~\ref{sec:alpha1} shows that the flow is harmless above $1$, where every law is a Poisson smoothing of the law at $1$, and that at $1$ itself the positivity is provable from the primes on windows as fine as the mean zero spacing. The obstacle is resolution and position at once: below the zero spacing the prime sum must cancel by more than a trivial bound allows, and below the line $\alpha = 1$ positivity at any fixed $\alpha$ is the quasi Riemann hypothesis.

The third route is the string. Under RH with simple zeros, $\mu$ is the hitting time of a Krein string with Dirichlet spectrum $\{\gamma^2\}$, by Kent's theorem. A construction of that string from the primes, that is, a speed measure on an interval whose spectral measure is the prime side of the explicit formula, would prove RH; this is a probabilistic form of the Hilbert Pólya programme. The obstacle is the usual one for Hilbert Pólya: nothing in the known structure produces the operator; the spectrum is known only through the zeros it is supposed to explain.

So the path, if there is one, runs through one of these three doors, and the two sanity tests say that no argument that fails to use the Euler product beyond the small $x$ regime can open any of them. What the present framework contributes is not a proof but a sharp statement of what is missing: an inequality that uses $\Lambda(n) \ge 0$ at every order of the ladder, failing for Davenport Heilbronn and for $H_t$ with $t < 0$, and reducing in the limit to the complete monotonicity of $A - P$. Finding one such inequality at a single order beyond the small $x$ regime, say at the first order where the prime term is no longer negligible, would be the first real step.

\section{Open problems}\label{sec:open}

The analysis points to a small number of precise questions.

\begin{enumerate}
\item[1.] \emph{The one sided bound.} Prove, for some fixed $\alpha < 1$ and all heights, that $-\Real\,\zeta'/\zeta(\alpha+i\tau) < \tfrac12\log(\tau/2\pi) + O(1/\tau)$, by a route that uses neither the convergent Euler product nor the location of zeros. By Theorem~\ref{thm:lagarias} this is the quasi Riemann hypothesis $\sigma_{\max} \le \alpha$; for a sequence $\alpha_n \downarrow \tfrac12$ it is RH, and by Theorem~\ref{thm:continuity} the sequence suffices.

\item[2.] \emph{A prime side proof at $\alpha = 1$.} The positivity of the Thorin density at $\alpha = 1$ is a theorem, but its proof here goes through the zeros. In the Mertens form \eqref{eq:mertensform} it says that the oscillating prime sum $\sum\Lambda(n)n^{-1}\cos(\tau\log n)$ is at most half the Mertens sum $\sum_{n\le\tau/2\pi}\Lambda(n)/n$, up to $\gamma_{\rm E}/2 + o(1)$. Find a proof from the primes, using no information about the zeros beyond the prime number theorem. Theorem~\ref{thm:fejer} proves it from Mertens' theorem, given the prime number theorem, for every smoothing of bandwidth up to $\log(\tau/2\pi)$, that is, on windows of one mean zero spacing; the problem is to pass to bandwidths growing faster than $\log\tau$, where the trivial bound fails by the factor $L/\log\tau$. Such a proof would be a model for the descent, since the same inequality one step to the left of $1$ is problem 1.

\item[3.] \emph{The negative mass.} Under RH the arithmetic factor $T_{\log N}$ in the tilted balance removes exactly the negative Thorin mass $-1.124$ of the archimedean law on $(\lambda, \lambda + 39.56)$ and adds nothing negative elsewhere. Quantify this directly: show that the Thorin side of a positive superposition of first passage times to the levels $\log n$, with weights $\Lambda(n)/(\sqrt n\log n)$, has negative part exactly matching \eqref{eq:thorinR} below $\gamma_1^2 + \lambda$.

\item[4.] \emph{The Cauchy flow.} The descent $\alpha \downarrow \tfrac12$ is a backward Cauchy flow on the zero counting measure, and $\sigma_{\max} - \tfrac12$ is its constant, with $\Lambda \le 2(\sigma_{\max}-\tfrac12)^2$ for the de Bruijn Newman constant. Is there a sharper relation between the two constants, and does the Cauchy flow admit an analogue of the Rodgers Tao argument~\cite{RodgersTao} that gave $\Lambda \ge 0$, which in Cauchy terms would be the statement that the flow cannot be continued below $\tfrac12$?

\item[5.] \emph{Below $\alpha = 1$ at finite height.} At any base point in the strip the law is within $1.5\times10^{-12}(\sqrt s + s)$ in log Laplace transform of an explicit GGC built from the zeros below $H$. Determine whether this proximity, combined with the heavy tail $\mathsf B_\alpha/\sqrt{\pi y}$ that every $\mu^{(\alpha)}$ with $\alpha > \tfrac12$ must have, constrains the possible Thorin measures of the remainder $R^{(\alpha)}$ beyond what the zero counting function gives.
\end{enumerate}

\section{Summary of the objects}

\begin{table}[htbp]
\centering
\small
\begin{tabular}{p{2.9cm}p{5.6cm}p{4.9cm}}
\toprule
Object & Unconditional & Under RH\\
\midrule
Law $\mu$ & self decomposable, absolutely continuous, unimodal; mean $0.0231$, s.d.\ $0.0061$ & GGC, law of $\sum E_\gamma/\gamma^2$\\
Canonical function $\Theta$ & $\sum_{\gamma>0}e^{(\rho-1/2)^2x}$; explicit formula \eqref{eq:weil}; level form \eqref{eq:levelform}; positive, nonincreasing & $\sum_{\gamma>0}e^{-\gamma^2x}$, completely monotone\\
Lévy measure & $\Theta(x)x^{-1}dx$, equation \eqref{eq:levydensity}: pole, stable, archimedean levels, minus first passages at $\log n$ & $\sum_\gamma e^{-\gamma^2x}x^{-1}dx$\\
Driver $D$ & jump density $-\Theta'(x)$, infinite activity & superposition of rate one compound Poisson with $\mathrm{Exp}(\gamma^2)$ jumps\\
Thorin measure $U$ & not a positive measure; on line part $\sum_{\rm on}\delta_{\gamma^2}$; boundary argument counts on line zeros & $\sum_{\gamma>0}\delta_{\gamma^2}$; distribution $N(\sqrt t) = \vartheta(\sqrt t)/\pi + 1 + S(\sqrt t)$\\
Smooth Thorin density & signed: $(4\pi\sqrt t)^{-1}(\Real\frac{\Gamma'}{\Gamma}(\frac14+\frac{i\sqrt t}2)-\log\pi)$, negative on $(0,39.56)$, mass $-1.124$, unit atom at $0$ & the same, plus $d[S(\sqrt t)]$\\
Tilted archimedean Thorin measure & exact: unit atom at $\lambda-\frac14$ plus the shifted smooth density, equation \eqref{eq:thorinR} & \\
Law at $\alpha = 1$, $1/(2\xi(1+\sqrt s))$ & GGC; Thorin measure absolutely continuous (prime number theorem), density $\pi^{-1}\Real\,\xi'/\xi(1+i\tau)$; tail $\mathsf B_1/\sqrt{\pi y}$, $\mathsf B_1 = 1+\frac{\gamma_{\rm E}}2-\frac12\log4\pi$ & Cauchy smoothing of the zero count with width $\frac12$\\
Prime side at $\alpha = 1$ & $U_\alpha = P_{\alpha-1}*U_1$ for $\alpha > 1$; Fejér smoothing of $U_1$ positive by Mertens for bandwidth $L \le \log(\tau/2\pi)+2\gamma_{\rm E}$ & positivity at all bandwidths\\
Base point law $\mu^{(\alpha)}$ & GGC for $\alpha\ge1$, Thorin measure a sum of Lorentzians \eqref{eq:lorentz}, distribution \eqref{eq:arctan}; heavy tail $\mathsf B_\alpha/\sqrt{\pi y}$ & GGC for every $\alpha \ge \frac12$; $U_\alpha \to \sum\delta_{\gamma^2}$ weakly as $\alpha\downarrow\frac12$\\
\bottomrule
\end{tabular}
\caption{The reciprocal Xi law and its measures.}
\end{table}

\appendix
\section{Numerical methods}\label{sec:numerics}

All computations use \texttt{mpmath} at 15 to 30 significant digits, with \texttt{numpy} for vectorised sums; the scripts are distributed with the paper. Zeros are computed with \texttt{zetazero}: the first $600$ for zero side sums, with the tail beyond the last zero $G$ replaced by $\int_G^\infty f(T)\,(2\pi)^{-1}\log(T/2\pi)\,dT$. Prime sums use exact von Mangoldt weights to $2\times10^5$. The explicit formula \eqref{eq:weil} is evaluated with the archimedean integral by adaptive quadrature and the prime sum truncated where $e^{-(\log n)^2/(4x)}$ is below machine precision; it agrees with the zero side to $3\times10^{-16}$ on $0.002 \le x \le 0.06$ (Figure~\ref{fig:law}). The cumulants are Taylor coefficients of $\log\xi(\tfrac12+w)$ from \texttt{taylor}, checked against zero sums to six digits.

The density of $\mu$ is computed by the inversion $f(y) = \pi^{-1}\int_0^\infty\Real[\varphi(-it)e^{-ity}]\,dt$ with $\varphi(-it) = \xi(\tfrac12)/\xi(\tfrac12+\sqrt{-it})$, trapezoidal rule with step $0.25$ on $[0,200]$, $2$ on $[200,3000]$ and $10$ on $[3000,30000]$; the transform is below $10^{-44}$ at the end. Base point products (Theorem~\ref{thm:product}) and Thorin integrals use the same zeros and tail. The values of $\Real\,\xi'/\xi$ on vertical lines use $\xi'/\xi(s) = 1/s + 1/(s-1) - \tfrac12\log\pi + \tfrac12\psi(s/2) + \zeta'/\zeta(s)$ with \texttt{zeta(s, derivative=1)}, on grids of step $0.05$ to $0.25$; minima are grid minima. Fejér weighted sums use the closed form weights $(1-\log n/L)_+$ and exact $\Lambda(n)$. Figures are produced by the same scripts with \texttt{matplotlib}.

\end{document}